\documentclass[11pt,a4paper]{article}
\usepackage[T1]{fontenc}
\usepackage{lmodern,amsmath,amssymb,amsthm,mathtools,microtype,array}
\usepackage[margin=27mm]{geometry}
\usepackage{enumitem,needspace,xcolor}
\usepackage{graphicx,tikz,float}
\usetikzlibrary{arrows.meta,calc,fit,positioning}
\usepackage{hyperref}
\hypersetup{colorlinks=true,linkcolor=blue!40!black,urlcolor=blue!40!black,
 citecolor=blue!40!black,pdftitle={Binary Triple Homogeneity Requires Exactly Three Binary Pair Stages}}
\newtheorem{theorem}{Theorem}[section]
\newtheorem{lemma}[theorem]{Lemma}
\newtheorem{corollary}[theorem]{Corollary}
\newtheorem{proposition}[theorem]{Proposition}
\newtheorem{claim}{Claim}
\theoremstyle{definition}\newtheorem{definition}[theorem]{Definition}

\theoremstyle{remark}\newtheorem{remark}[theorem]{Remark}
\newenvironment{claimproof}
 {%
  \begin{proof}[Proof of the claim]}
 {\end{proof}}
\numberwithin{equation}{section}
\newcommand{\Hom}{\mathbf{Hom}}
\newcommand{\BT}{\mathrm{BT}}
\newcommand{\Ccal}{\mathcal C}
\newcommand{\N}{\mathbb N}
\newcommand{\NatLt}[1]{\mathbb N_{<#1}}
\newcommand{\Bits}{\{0,1\}}
\newcommand{\Ninf}{\mathbb N_\infty}
\makeatletter
\newcommand{\X}{\@ifnextchar^{\X@power}{[\N]^\N}}
\def\X@power^#1{\bigl([\N]^\N\bigr)^{#1}}
\makeatother
\newcommand{\dom}{\operatorname{dom}}
\newcommand{\restr}{\mathbin{\upharpoonright}}
\newcommand{\abs}[1]{\lvert#1\rvert}
\setlist[enumerate]{itemsep=3pt,topsep=5pt}
\allowdisplaybreaks[1]
\title{Borel--Tukey morphisms between Ramsey homogeneity relations}
\author{Tatsuya Goto\\[2pt]
\small Institute of Discrete Mathematics and Geometry, TU Wien,\\
\small Wiedner Hauptstra\ss{}e 8--10/104, 1040 Wien, Austria\\
\small e-mail: \texttt{goto.tatsuya@icloud.com}}
\date{\today}
\begin{document}
\maketitle
\vspace{-2em}
\begin{abstract}

We study Borel--Tukey morphisms between almost-homogeneity relations.
For finite \(r,s\ge1\) and \(m,n\ge2\), we give the complete
one-stage classification: \(\Hom^r_m\to_{\BT}\Hom^s_n\) holds exactly
when \(r>s\), or when \(r=s\) and \(m\ge n\).
Our main result is that binary triple homogeneity requires exactly
three sequential binary pair stages:
\(
 \Hom^2_2;\Hom^2_2;\Hom^2_2\longrightarrow_{\BT}\Hom^3_2\) and \(
 \Hom^2_k;\Hom^2_2\not\longrightarrow_{\BT}\Hom^3_2
\)
for every finite positive integer \(k\).
We also construct Borel morphisms
\(\Hom^2_2;\Hom^2_3\longrightarrow_{\BT}\Hom^3_2\) and
\(\Hom^2_2;\Hom^2_2\longrightarrow_{\BT}\Hom^2_k\) for every finite
positive \(k\).
\end{abstract}

\tableofcontents

\section{Introduction}\label{sec:introduction}
The infinite Ramsey theorem says that every finite coloring of the
\(r\)-element subsets of \(\N\) has an infinite homogeneous set.
We study the relations obtained by taking such colorings as challenges
and infinite almost homogeneous sets as responses, allowing finitely
many initial exceptions. Although Ramsey's theorem guarantees that
responses exist, Borel--Tukey morphisms compare how uniformly one such
homogeneity problem can be transformed into another.

Borel--Tukey morphisms can distinguish relations whose cardinal
invariants are equal. For a finite \(m\ge2\), let
\[
 \mathbf R_m=((\NatLt{m})^{\N},[\N]^\N,R_m),
 \qquad
 cR_mA\quad\Longleftrightarrow\quad
 c\restr A\text{ is almost constant}.
\]
For finite \(m,n\ge2\),
\[
 \mathbf R_m\longrightarrow_{\BT}\mathbf R_n
 \quad\Longleftrightarrow\quad n\le m.
\]
When \(n\le m\), the morphism is immediate; Mildenberger proved that
no Borel morphism exists when \(m<n\)~\cite{Mildenberger}.

The higher-dimensional homogeneity relations considered here cannot be
separated by their cardinal invariants. Blass showed that, for all finite
\(r,m\ge2\), the least size of a family of infinite sets that contains an
almost homogeneous set for every \(m\)-coloring of \([\N]^r\) is
\[
 \max\{\mathfrak r_\sigma,\mathfrak d\},
\]
where \(\mathfrak r_\sigma\) and \(\mathfrak d\) are the
\(\sigma\)-reaping and dominating numbers~\cite{Blass}. Thus cardinal
invariants cannot distinguish the relations \(\Hom^r_m\). This motivates
studying the morphisms themselves. Our main result determines the exact
number of binary pair stages needed for binary triple homogeneity. We
also give the complete classification of one-stage morphisms between
the finite-arity relations \(\Hom^r_m\). To formulate these results
precisely, we identify subsets of a countable set with their
characteristic functions in the product topology. We use the convention
that \(0\in\N\). For \(n\in\N\), write
\[
 \NatLt{n}=\{k\in\N:k<n\}.
\]
For \(A\subseteq\N\) and \(n\in\N\), use the notation
\[
 \begin{aligned}
 [A]^n&=\{s\subseteq A:\abs{s}=n\},\\
 [A]^{<\infty}&=\{s\subseteq A:s\text{ is finite}\},\\
 [A]^\N&=\{B\subseteq A:B\text{ is infinite}\}.
 \end{aligned}
\]
For \(r,m\ge1\), put
\[
 \Ccal^r_m=(\NatLt{m})^{[\N]^r} = \{ c : c\text{ is a function from }[\N]^r\text{ to }\NatLt{m} \}.
\]
For \(c\in\Ccal^r_m\) and \(A\in\X\), write \(cH^r_m A\) if \(A\) is \emph{almost homogeneous} for \(c\): there are \(N\in\N\) and \(\varepsilon\in\NatLt{m}\) with
\[
 c(s)=\varepsilon\qquad(s\in[A\cap[N,\infty)]^r).
\]
Thus vertices below the cutoff \(N\) may be removed. A set is
\emph{exactly homogeneous} if none
are removed. Define
\[
 \Hom^r_m=(\Ccal^r_m,\X,H^r_m).
\]
The superscript is the arity and the subscript is the number of colors.
For \(A\in\X\), let \(e_A:\N\to A\) be its increasing
enumeration, starting at zero. For
\(s, B \subseteq \N\), write \(s<B\) if every element of \(s\) is below every element of \(B\); this holds vacuously when \(s=\varnothing\).

A morphism \(\Hom^r_m\to\Hom^s_n\) consists of maps
\[
 \alpha:\Ccal^s_n\to\Ccal^r_m,
 \qquad\beta:\X\to\X
\]
satisfying
\begin{equation}\label{eq:one-stage}
 \alpha(c)H^r_m A\quad\Longrightarrow\quad cH^s_n\beta(A).
\end{equation}
The response \(\beta(A)\) cannot use \(c\). The morphism is
Borel--Tukey if both maps are Borel.

For relations \(\mathbf A=(A_-,A_+,R)\) and
\(\mathbf B=(B_-,B_+,S)\), sequential composition has challenges
\((a,f)\in A_-\times{}^{A_+}B_-\), responses
\((x,y)\in A_+\times B_+\), and relation
\[
 (a,f)\,T\,(x,y)\quad\Longleftrightarrow\quad aRx\ \land\ f(x)Sy.
\]
To avoid putting a topology on an unrestricted function space, we
state Borel regularity directly for the stage maps. A
morphism from $\Hom^2_2; \Hom^2_2; \Hom^2_2$ to \(\Hom^3_2\) consists of jointly Borel maps
\begin{align*}
 \alpha_1&:\Ccal^3_2\to\Ccal^2_2,\\
 \alpha_2&:\Ccal^3_2\times\X\to\Ccal^2_2,\\
 \alpha_3&:\Ccal^3_2\times\X^2\to\Ccal^2_2,\\
 \Gamma&:\X^3\to\X
\end{align*}
with
\begin{equation}\label{eq:three-stage-condition}
 \begin{split}
 &\alpha_1(c)H^2_2 A\ \land\ \alpha_2(c,A)H^2_2 B
       \ \land\ \alpha_3(c,A,B)H^2_2 C\\
 &\hspace{35mm}\Longrightarrow\quad cH^3_2\Gamma(A,B,C).
 \end{split}
\end{equation}

Similarly, a two-stage morphism
\(\Hom^2_k;\Hom^2_\ell\to_{\BT}\Hom^s_n\) is specified by
jointly Borel maps
\[
 \alpha_1:\Ccal^s_n\to\Ccal^2_k,\qquad
 \alpha_2:\Ccal^s_n\times\X\to\Ccal^2_\ell,\qquad
 \Gamma:\X^2\to\X
\]
such that
\begin{equation}\label{eq:two-stage-condition}
 \alpha_1(c)H^2_k A\ \land\ \alpha_2(c,A)H^2_\ell B
 \quad\Longrightarrow\quad cH^s_n\Gamma(A,B).
\end{equation}

The main results are as follows.
\begin{theorem}\label{thm:overview}
\begin{enumerate}
\item (Theorem~\ref{thm:complete-classification}) For finite
\(r,s\ge1\) and \(m,n\ge2\),
\[
 \Hom^r_m\longrightarrow_{\BT}\Hom^s_n
 \quad\Longleftrightarrow\quad
 r>s\ \text{ or }\ (r=s\text{ and }m\ge n).
\]
\item (Theorem~\ref{thm:three-stage})
\[
 \Hom^2_2;\Hom^2_2;\Hom^2_2\longrightarrow_{\BT}\Hom^3_2.
\]
\item (Theorem~\ref{thm:two-stage-triples})
\[
 \Hom^2_2;\Hom^2_3\longrightarrow_{\BT}\Hom^3_2.
\]
\item (Theorem~\ref{thm:two-stage-palettes}) For every finite positive integer \(k\),
\[
 \Hom^2_2;\Hom^2_2\longrightarrow_{\BT}\Hom^2_k.
\]
\item (Theorem~\ref{thm:two-stage-negative}) For every finite positive integer \(k\),
\[
\Hom^2_k;\Hom^2_2 \not\longrightarrow_\BT\Hom^3_2.
\]
\end{enumerate}
\end{theorem}

The ultimate goal of this line of work is to characterize completely
the finite parameter sequences
\[
 ((r_1,m_1),\ldots,(r_\ell,m_\ell),(s,n))
\]
for which
\[
 \Hom^{r_1}_{m_1};\cdots;\Hom^{r_\ell}_{m_\ell}
 \longrightarrow_{\BT}\Hom^s_n
\]
holds.

\section{Normalization lemmas}\label{sec:normalization}

We develop normalization tools that replace Borel responses by
continuous online selectors and make Borel queries continuous after
reindexing.

We use the following form of the Galvin--Prikry theorem: for every
infinite \(R\subseteq\N\) and every Borel map
\(f:[R]^\N\to\Bits\), there is an infinite \(R'\subseteq R\) on
which \(f\) is constant, meaning that \(f\restr[R']^\N\) is
constant \cite{GP}. Shrinking \(R\) successively handles finitely many
such requirements.

\begin{lemma}\label{lem:compact}
Let \(\beta:[\N]^\N\to\Bits^\N\) be Borel. There are an
infinite \(M\subseteq\N\) and a continuous map
\[
  F:\Bits^M\longrightarrow\Bits^\N
\]
such that \(F(A)=\beta(A)\) for every \(A\in[M]^\N\).
\end{lemma}

\begin{proof}
Put \(R_{-1}=\N\). Recursively choose integers
\(m_0<m_1<\cdots\) and infinite sets \(R_n\) as follows.
At stage \(n\), choose \(m_n\in R_{n-1}\) and set
\[
  S_n=\{m_0,\ldots,m_n\}.
\]
Start with the infinite tail
\(R_{n-1}\cap(m_n,\infty)\). For each of the finitely many
\(t\subseteq S_n\), apply Galvin--Prikry to the Borel map
\[
  B\longmapsto
  \begin{cases}
   1,& n\in\beta(t\cup B),\\
   0,& n\notin\beta(t\cup B).
  \end{cases}
\]
After finitely many refinements, obtain an infinite
\(R_n\subseteq R_{n-1}\cap(m_n,\infty)\) and values
\(\varepsilon_n(t)\in\Bits\), for \(t\subseteq S_n\), such that
\begin{equation}\label{eq:bitdecision}
  \varepsilon_n(t)=
  \begin{cases}
   1,& n\in\beta(t\cup B),\\
   0,& n\notin\beta(t\cup B)
  \end{cases}
  \qquad(t\subseteq S_n,\ B\in[R_n]^\N).
\end{equation}
Earlier requirements remain true under subsequent refinements.

Let \(M=\{m_n:n\in\N\}\). By construction,
\(M\setminus S_n\subseteq R_n\). For every subset \(D\subseteq M\),
including finite \(D\), define
\begin{equation}\label{eq:extension}
  F(D)(n)=\varepsilon_n(D\cap S_n).
\end{equation}
The \(n\)-th coordinate depends only on the finite set of input
coordinates \(S_n\). Hence each coordinate map is continuous, and so
is \(F\).

If \(A\in[M]^\N\), the set \(A\setminus S_n\) is infinite and
is contained in \(R_n\). Applying \eqref{eq:bitdecision} with
\(t=A\cap S_n\) and \(B=A\setminus S_n\) gives
\[
  F(A)(n)=\varepsilon_n(A\cap S_n)
  =\begin{cases}
    1,& n\in\beta(A),\\
    0,& n\notin\beta(A).
   \end{cases}
\]
This holds for every \(n\), proving the assertion.
\end{proof}

Let \(\beta\), \(F\), and \(M\) be as in
Lemma~\ref{lem:compact}.
For finite \(t\subseteq M\), write \(B_t=F(t)\). We use two
consequences of continuity. First, if \(E\subseteq\N\) is finite,
then, for all sufficiently
large \(q\),
\begin{equation}\label{eq:uniformtail}
  F(t\cup C)\cap E=B_t\cap E
  \qquad\bigl(C\subseteq M\cap[q,\infty)\bigr).
\end{equation}
Indeed, continuity of the map \(D\mapsto F(D)\cap E\) at \(t\)
provides a finite set \(K\subseteq M\) such that
\[
  D\cap K=t\cap K\quad\Longrightarrow\quad
  F(D)\cap E=F(t)\cap E.
\]
Choose \(q\) larger than every element of \(K\cup t\). Then every
\(C\subseteq M\cap[q,\infty)\) is disjoint from \(K\cup t\), so
\((t\cup C)\cap K=t\cap K\), which proves \eqref{eq:uniformtail}.
The statement holds for finite and infinite \(C\), as well as
\(C=\varnothing\). Second, if \(t_j\) are the increasing finite
initial segments of \(A\in[M]^\N\), then
	\begin{equation}\label{eq:branchlimit}
  B_{t_j}\longrightarrow F(A)=\beta(A)
  \quad\text{in }\Bits^\N\text{ as }j\to\infty,
\end{equation}
because \(t_j\to A\) in \(\Bits^M\).

The same tail-agreement argument works for every continuous map on
\(\Bits^M\). Lemma~\ref{lem:trim} below thins the outputs so that finite
inputs have finite outputs, without changing the needed homogeneity
implications.

\begin{lemma}\label{lem:trim}
Assume that \(\beta:\X\to\X\) is a Borel map, and take \(M,F\) as in Lemma~\ref{lem:compact}. Define
\begin{equation}\label{eq:trim}
 F^*(D)=\{b\in F(D):\abs{F(D)\cap\NatLt{b}}\in D\}
 \qquad(D\subseteq M).
\end{equation}
Then \(F^*\) is continuous, \(F^*(\varnothing)=\varnothing\), and \(\abs{F^*(t)}\le\abs t\) for finite \(t\). For infinite \(A\subseteq M\),
\[
 F^*(A)=e_{\beta(A)}[A]\in[\beta(A)]^\N.
\]
\end{lemma}
\begin{proof}
Membership of \(b\) in \eqref{eq:trim} depends only on
\(F(D)\restr\NatLt{b+1}\) and \(D\restr\NatLt{b+1}\), so \(F^*\) is continuous.
The map \(b\mapsto\abs{F(D)\cap\NatLt{b}}\) is injective on \(F(D)\).
Thus each element of a finite \(t\) selects at most one element of
\(F(t)\), giving \(\abs{F^*(t)}\le\abs t\). For infinite \(A\), write
\(\beta(A)=\{b_0<b_1<\cdots\}\). Since
\(\abs{\beta(A)\cap\NatLt{b_n}}=n\), definition~\eqref{eq:trim} gives
\[
 F^*(A)=\{b_n:n\in A\}=e_{\beta(A)}[A]. \qedhere
\]
\end{proof}

\begin{definition}
An \emph{online selector} is a map \(\gamma:\Bits^\N\to\Bits^\N\) such that
\(\gamma(D)\subseteq D\) for every \(D\subseteq\N\),
\(\gamma(D)\) is infinite whenever \(D\) is infinite, and
\begin{equation}\label{eq:online}
 D\cap\NatLt{k+1}=E\cap\NatLt{k+1}
 \quad\Longrightarrow\quad
 \gamma(D)\cap\NatLt{k+1}=\gamma(E)\cap\NatLt{k+1}.
\end{equation}
Thus retaining \(k\) depends only on input coordinates at most \(k\),
so the map is continuous.
\end{definition}

\begin{lemma}[Sparse block normalization]\label{lem:block-normal}
For every Borel \(\beta:\X\to\X\), there are a strictly increasing
map \(u:\N\to\N\) with \(u(n+1)>u(n)+1\), a continuous map
\(F:\Bits^\N\to\Bits^\N\), and pairwise disjoint finite sets
\((Q_n)_{n\in\N}\), with \(Q_n\subseteq[n,\infty)\), such that:
\begin{enumerate}
\item \(F(C)\in[\beta(u[C])]^\N\) for every \(C\in\X\);
\item \(F(D)\subseteq\bigcup_{n\in D}Q_n\) for every
  \(D\subseteq\N\);
\item \(F(D)\cap Q_n=F(D\cap\NatLt{n+1})\cap Q_n\) for every
  \(D\subseteq\N\) and \(n\in\N\).
\end{enumerate}
\end{lemma}
\begin{proof}
Apply Lemmas~\ref{lem:compact} and \ref{lem:trim} to obtain an infinite
\(M\subseteq\N\) and a continuous map \(G:\Bits^M\to\Bits^\N\)
such that \(G(\varnothing)=\varnothing\), \(G(t)\) is finite for
finite \(t\), and
\[
 G(A)\in[\beta(A)]^\N
 \qquad(A\in[M]^\N).
\]
For finite \(t\subseteq M\), write \(B_t=G(t)\).

Recursively construct \(X=\{x_n:n\in\N\}\subseteq M\) and the
sets \(Q_n\). At stage \(n\), put
\[
 X_n=\{x_i:i\in\NatLt{n}\},\qquad
 D_n=\bigcup_{i\in\NatLt{n}}Q_i,
\]
and set
\[
 E_n=D_n\cup\NatLt{n}\cup\bigcup_{t\subseteq X_n}B_t.
\]
Continuity at the finitely many \(t\subseteq X_n\) allows choosing
\(x_n\in M\) sufficiently large that
\begin{equation}\label{eq:preserve}
 G(t\cup C)\cap E_n=B_t\cap E_n
 \quad\left(t\subseteq X_n,
 C\subseteq M\cap[x_n,\infty)\right).
\end{equation}
Also require \(x_n>x_{n-1}+1\) for \(n>0\). 

Since \(B_t\subseteq E_n\), equation \eqref{eq:preserve}, applied with
\(C=\{x_n\}\), gives
\[
 B_{t\cup\{x_n\}}\cap E_n=B_t.
\]
In particular, \(B_t\subseteq B_{t\cup\{x_n\}}\) and
\((B_{t\cup\{x_n\}}\setminus B_t)\cap E_n=\varnothing\). Therefore
\[
 Q_n=\bigcup_{t\subseteq X_n}
          (B_{t\cup\{x_n\}}\setminus B_t)
\]
is finite and disjoint from \(E_n\). By definition,
\(D_n=\bigcup_{i\in\NatLt{n}}Q_i\), and \(D_n\subseteq E_n\).
Consequently, \(Q_n\cap Q_i=\varnothing\) for every \(i<n\); Thus, the sets \((Q_n)_{n\in\N}\) are pairwise
disjoint. Also \(\NatLt{n}\subseteq E_n\), so
\(Q_n\cap\NatLt{n}=\varnothing\). Hence \(Q_n\subseteq[n,\infty)\).

Let \(u(n)=x_n\) and define \(F(D)=G(u[D])\). The map \(F\) is
continuous.

We now prove the item (2) of this lemma.
Assume $D$ is finite and let \(D=\{n_0<\cdots<n_{r-1}\}\) and put
\[
 t_\ell=\{x_{n_j}:j\in\NatLt{\ell}\}
 \qquad(0\le\ell\le r).
\]
Then \(t_0=\varnothing\), \(t_r=u[D]\), and
\(t_\ell\subseteq X_{n_\ell}\) for \(\ell<r\). Hence the definition of
\(Q_{n_\ell}\) gives
\[
 B_{t_{\ell+1}}\setminus B_{t_\ell}\subseteq Q_{n_\ell}
 \qquad(\ell<r).
\]
The sets \(B_{t_\ell}\) are increasing, and
\(B_{t_0}=B_\varnothing=\varnothing\). Therefore
\[
 \begin{aligned}
 F(D)=B_{t_r}
 &=\bigcup_{\ell\in\NatLt{r}}
       (B_{t_{\ell+1}}\setminus B_{t_\ell})\\
 &\subseteq\bigcup_{\ell\in\NatLt{r}}Q_{n_\ell}
  =\bigcup_{n\in D}Q_n.
 \end{aligned}
\]
For infinite \(D\), continuity and the same monotonicity give
\(F(D)=\bigcup_jF(D\cap\NatLt{j})\), so the inclusion holds for every
\(D\subseteq\N\).

If \(C\) is infinite, the defining property of
\(G\) gives \(F(C)\in[\beta(u[C])]^\N\).

Finally, fix \(D\subseteq\N\) and \(n\in\N\), and put
\[
 t=u[D\cap\NatLt{n+1}],\qquad
 C=u[D\setminus\NatLt{n+1}].
\]
Then \(t\subseteq X_{n+1}\),
\(C\subseteq M\cap[x_{n+1},\infty)\), and \(u[D]=t\cup C\).
Moreover, \(Q_n\subseteq D_{n+1}\subseteq E_{n+1}\). Therefore
\eqref{eq:preserve}, applied at stage \(n+1\), gives
\[
 \begin{aligned}
 F(D)\cap Q_n
 &=G(t\cup C)\cap Q_n\\
 &=G(t)\cap Q_n\\
 &=F(D\cap\NatLt{n+1})\cap Q_n.
 \end{aligned}
\]
\end{proof}

\begin{corollary}[Online normalization]\label{lem:normal}
For a Borel \(\beta:\X\to\X\), there exist an infinite coinfinite
\(X\subseteq\N\), a continuous \(F^*:\Bits^X\to\Bits^\N\), and a
finite-to-one map \(h:\N\to\N\) such that:
\begin{enumerate}
\item \(F^*(A)\in[\beta(A)]^\N\) for every infinite \(A\subseteq X\);
\item \(h[F^*(D)]\subseteq D\) for every \(D\subseteq X\);
\item writing \(e=e_X\), the map
\begin{equation}\label{eq:gamma}
 \gamma(D)=e^{-1}[h[F^*(e[D])]]\qquad(D\subseteq\N)
\end{equation}
is an online selector.
\end{enumerate}
\end{corollary}
\begin{proof}
Take \(u,F,(Q_n)_n\) from Lemma~\ref{lem:block-normal}, put
\(X=u[\N]\), and let \(e=e_X=u\). Since
\(u(n+1)>u(n)+1\), we have \(u(n)+1\notin X\) for every \(n\in\N\).
So \(\N\setminus X\) is infinite. Define
\[
 F^*(A)=F(e^{-1}[A])\qquad(A\subseteq X).
\]
Define \(h(b)=e(n)\) for \(b\in Q_n\), and map the remaining points
injectively into \(\N\setminus X\). Then \(h\) is finite-to-one and
\(h^{-1}(\{e(n)\})=Q_n\). The first two conclusions follow from the
first two clauses of Lemma~\ref{lem:block-normal}.

For \(D\subseteq\N\), we have \(F^*(e[D])=F(D)\), and the definition
of \(h\) gives
\[
 n\in\gamma(D)\quad\Longleftrightarrow\quad F(D)\cap Q_n\ne\varnothing.
\]
Using the block-support property in the second clause of
Lemma~\ref{lem:block-normal} and the fact that the blocks \(Q_k\) are pairwise disjoint, it follows that
\(\gamma(D)\subseteq D\).
The last clause of Lemma~\ref{lem:block-normal} gives the online property. If \(D\) is
infinite, then \(F(D)\) is infinite. Also, since \(h\) is finite-to-one,
\(\gamma(D)\) is infinite.
\end{proof}

\begin{corollary}[Tail-controlled online normalization]\label{lem:tail-normal}
For every Borel \(\beta:\X\to\X\), there are a strictly increasing
map \(u:\N\to\N\), a continuous map
\(F:\Bits^\N\to\Bits^\N\), and a sequence \((Q_n)_{n\in\N}\) of
pairwise disjoint finite sets with \(Q_n\subseteq[n,\infty)\) such that:
\begin{enumerate}
\item \(F(C)\in[\beta(u[C])]^\N\) for every \(C\in\X\);
\item \(F(D)\subseteq\bigcup_{n\in D}Q_n\) for every
  \(D\subseteq\N\);
\item \(F(D)\cap Q_n=F(D\cap\NatLt{n+1})\cap Q_n\) for every
  \(D\subseteq\N\) and \(n\in\N\).
\end{enumerate}
\end{corollary}
\begin{proof}
This is Lemma~\ref{lem:block-normal} with the gap condition on \(u\)
discarded.
\end{proof}

\begin{corollary}[Overlap normalization]\label{lem:overlap}
For every Borel \(\beta:\X\to\X\), there are an infinite \(X\subseteq\N\) and a finite-to-one \(h:\N\to\N\) such that \(\beta\restr[X]^\N\) is continuous and
\[
 \abs{h[\beta(A)]\cap A}=\infty
 \qquad(A\in[X]^\N).
\]
\end{corollary}
\begin{proof}
Use the \(X,F^*,h\) produced in the proof of
Corollary~\ref{lem:normal}. Fix an infinite \(A\subseteq X\). Then
\(F^*(A)\) is an infinite subset of \(\beta(A)\). Since \(h\) is
finite-to-one, \(h[F^*(A)]\) is infinite: otherwise
\[
 F^*(A)\subseteq\bigcup_{m\in h[F^*(A)]}h^{-1}(\{m\})
\]
would express \(F^*(A)\) as a subset of a finite union of finite sets.
Moreover, the second conclusion of Corollary~\ref{lem:normal} gives
\(h[F^*(A)]\subseteq A\), while \(F^*(A)\subseteq\beta(A)\) gives
\(h[F^*(A)]\subseteq h[\beta(A)]\). Consequently,
\[
 h[F^*(A)]\subseteq h[\beta(A)]\cap A,
\]
and the set on the left is infinite. This proves the overlap assertion.
In the construction of Lemma~\ref{lem:block-normal}, the set
\(X=u[\N]\) is contained in the set \(M\) supplied by
Lemma~\ref{lem:compact}. Hence \(\beta\restr[X]^\N\) agrees with the
continuous extension obtained there.
\end{proof}
This does not say that \(\beta(A)\subseteq X\); instead, it gives
infinite overlap after a finite-to-one pullback.

\begin{lemma}[A full coloring cube in a comeager set]\label{lem:free-cube}
Let \(r,k\ge1\). If \(D\subseteq\Ccal^r_k\) is comeager, there are an
infinite \(Y=\{y_0<y_1<\cdots\} \subseteq \N\) and a continuous map
\(\Psi:\Ccal^r_k\to D\) such that
\begin{equation}\label{eq:free-cube}
 \Psi(c)(y_{i_0},\ldots,y_{i_{r-1}})
   =c(i_0,\ldots,i_{r-1})
 \qquad(i_0<\cdots<i_{r-1}).
\end{equation}
Every coordinate outside \([Y]^r\) is constant as \(c\) varies.
\end{lemma}
\begin{proof}
Put \(I=[\N]^r\) and \(K=\NatLt{k}\). For a finite partial map
\(p:I\rightharpoonup K\), let
\[
 \mathcal N(p)=\{c\in K^I:p\subseteq c\}.
\]
Choose dense open sets \(U_n\subseteq K^I\) such that
\[
 \bigcap_{n\in\N}U_n\subseteq D.
\]

We recursively construct
\[
 Y_n=\{y_i:i\in\NatLt{n}\},\qquad
 p_n:I\rightharpoonup K,
\]
where \(p_n\) is finite, \(p_n\subseteq p_{n+1}\), and
\begin{equation}\label{eq:free-cube-stage}
 \dom(p_n)\cap[Y_n]^r=\varnothing,
 \qquad
 \mathcal N(p_{n+1}\cup\sigma)\subseteq U_n
 \quad\bigl(\sigma\in K^{[Y_{n+1}]^r}\bigr).
\end{equation}
Start with \(Y_0=\varnothing\) and \(p_0=\varnothing\). Given
\(Y_n,p_n\), choose
\[
 y_n>\max\left(Y_n\cup\bigcup\dom(p_n)\right),
 \qquad Y_{n+1}=Y_n\cup\{y_n\},
\]
where \(\max\varnothing=-1\). Then
\(\dom(p_n)\cap[Y_{n+1}]^r=\varnothing\).

Enumerate the finite set \(K^{[Y_{n+1}]^r}\) as
\(\{\sigma_j:j\in\NatLt{N}\}\). Put \(p^{(0)}=p_n\). Given
\(p^{(j)}\), density and openness of \(U_n\) give a finite partial map
\(q_j\) such that
\[
 p^{(j)}\cup\sigma_j\subseteq q_j,
 \qquad \mathcal N(q_j)\subseteq U_n.
\]
Set
\[
 p^{(j+1)}=p^{(j)}\cup
 \bigl(q_j\restr(I\setminus[Y_{n+1}]^r)\bigr),
 \qquad p_{n+1}=p^{(N)}.
\]
For every \(j\in\NatLt{N}\), the map
\(p_{n+1}\cup\sigma_j\) extends \(q_j\); hence
\eqref{eq:free-cube-stage} holds.

Let
\[
 Y=\{y_n:n\in\N\},\qquad p_\infty=\bigcup_{n\in\N}p_n.
\]
If \(s\in\dom(p_{n+1})\setminus[Y_{n+1}]^r\), then
choose \(z\in s\setminus Y_{n+1}\). At stage \(n+1\), the point
\(y_{n+1}\) is chosen larger than every member of every set in
\(\dom(p_{n+1})\), so \(z<y_{n+1}\). Since \(z\notin Y_{n+1}\) and
\(y_m\ge y_{n+1}>z\) for every \(m\ge n+1\), we have \(z\notin Y\),
and hence \(s\notin[Y]^r\). Thus
\(\dom(p_\infty)\cap[Y]^r=\varnothing\). Extend \(p_\infty\) to a map
\[
 a:I\setminus[Y]^r\longrightarrow K.
\]
For \(c\in K^I\), define
\[
 \Psi(c)(s)=
 \begin{cases}
  c(e_Y^{-1}[s]),&s\in[Y]^r,\\
  a(s),&s\in I\setminus[Y]^r.
 \end{cases}
\]
Each output coordinate is either constant or an input coordinate, so
\(\Psi\) is continuous and satisfies \eqref{eq:free-cube}. For every
\(n\), the coloring \(\Psi(c)\) extends
\(p_{n+1}\cup(\Psi(c)\restr[Y_{n+1}]^r)\); hence
\[
 \Psi(c)\in\bigcap_{n\in\N}U_n\subseteq D.
\]
\end{proof}

\section{Complete one-stage classification}\label{sec:complete-classification}

We first classify one-stage morphisms at all finite arities. The
classification is lexicographic: arity takes priority, and at equal
arity the source must have at least as many colors as the target.

\begin{theorem}\label{thm:complete-classification}
For finite integers \(r,s\ge1\) and \(m,n\ge2\),
\begin{equation}\label{eq:complete-classification}
 \Hom^r_m\longrightarrow_{\BT}\Hom^s_n
 \quad\Longleftrightarrow\quad
 r>s\ \ \text{or}\ \ (r=s\text{ and }m\ge n).
\end{equation}
\end{theorem}

We divide the proof into the positive reductions, a normalization
lemma, and two nonreducibility arguments for the converse: colors cannot
be increased at fixed arity, and arity cannot be increased.

\subsection{The positive reductions}

At equal arity, \(n\le m\) gives the identity response: simply regard
an \(n\)-coloring as an \(m\)-coloring.  The reduction which lowers
arity is the following one.  Given \(c:[\N]^s\to\NatLt{n}\), define
\begin{equation}\label{eq:arity-lowering-query}
 q_c(u)=
 \begin{cases}
 0,&c\restr[u]^s\text{ is constant},\\
 1,&\text{otherwise},
 \end{cases}
 \qquad (u\in[\N]^{s+1}).
\end{equation}
Use the identity response.  No infinite set can be homogeneous for
\(q_c\) in color \(1\): an infinite \(c\)-homogeneous subset, supplied
by Ramsey's theorem, would give an \((s+1)\)-set of color \(0\).
If a set is \(q_c\)-homogeneous in color \(0\), any two of its
\(s\)-subsets can be joined by a chain which changes one point at a
time.  Consecutive members of the chain lie in a common
\((s+1)\)-set, so they have the same \(c\)-color.  Hence the whole set
is \(c\)-homogeneous.  Applying the same argument to a tail handles
almost homogeneity.  Thus
\begin{equation}\label{eq:one-step-arity-lowering}
 \Hom^{s+1}_2\longrightarrow_{\BT}\Hom^s_n.
\end{equation}
The challenge map is continuous.  For \(r>s+1\), assign an \(r\)-set
the \(q_c\)-color of its first \(s+1\) elements.  Every
\((s+1)\)-subset of an infinite set can be extended by larger points,
so this gives the required reduction for every \(r>s\) and \(m\ge2\).

\subsection{All-arity continuous normalization}

The negative direction uses the following consequence of Online
normalization and the full-cube argument.

\begin{lemma}\label{lem:all-arity-normalization}
If
\(\Hom^r_m\to_{\BT}\Hom^s_n\), then there is such a morphism
\((p,\gamma)\) in which both maps are continuous and
\begin{equation}\label{eq:all-arity-selector}
 \gamma(A)\subseteq A\qquad(A\in\X).
\end{equation}
Moreover, \(\gamma\) extends to an online selector on \(\Bits^\N\).
\end{lemma}
\begin{proof}
Start with a Borel morphism \((\alpha,\beta)\), and apply
Lemma~\ref{lem:normal} to \(\beta\).  Retain the resulting
\(X,F^*,h\), put \(e=e_X\), and let \(\gamma\) be the selector from
\eqref{eq:gamma}.  For a target coloring \(c\in\Ccal^s_n\), define a
coloring \(\bar c\in\Ccal^s_n\) on the original output space by
\begin{equation}\label{eq:all-arity-pullback}
 \bar c(v)=
 \begin{cases}
 c\bigl(e^{-1}[h[v]]\bigr),
    &h[v]\in[X]^s,\\
 0,&\text{otherwise}.
 \end{cases}
\end{equation}
Pull \(\alpha(\bar c)\) back along \(e\).  If \(A\) is almost
homogeneous for this new query, then \(e[A]\) is almost homogeneous
for \(\alpha(\bar c)\), and hence \(\beta(e[A])\) is almost
homogeneous for \(\bar c\).  Its infinite subset \(F^*(e[A])\) has
the same property.  Since \(h\) is finite-to-one, discarding finitely
many further points shows that
\[
 e^{-1}[h[F^*(e[A])]]=\gamma(A)
\]
is almost homogeneous for \(c\).  We have therefore obtained a Borel
challenge map and the continuous online selector \(\gamma\).

It remains to make the challenge map continuous.  Let \(D\) be a
comeager \(G_\delta\) set on which it is continuous.  Apply
Lemma~\ref{lem:free-cube} with arity \(s\) and \(n\) colors to obtain an
infinite
\(Y=\{y_0<y_1<\cdots\}\) and a continuous map
\begin{equation}\label{eq:general-free-cube}
 \Psi:\Ccal^s_n\longrightarrow D,
 \qquad
 \Psi(c)(y_{i_0},\ldots,y_{i_{s-1}})
   =c(i_0,\ldots,i_{s-1}).
\end{equation}
Pull the query back along \(e_Y\), and replace the response by
\begin{equation}\label{eq:all-arity-final-response}
 A\longmapsto e_Y^{-1}[\gamma(e_Y[A])].
\end{equation}
The selector property puts the old output inside \(Y\), so
\eqref{eq:general-free-cube} preserves correctness.  Both maps are
continuous, and the response is again an infinitude-preserving online
selector.
\end{proof}

\subsection{Colors cannot be increased at fixed arity}

Suppose first that \(r\ge2\), \(n>m\ge2\), and that
\begin{equation}\label{eq:forbidden-fixed-arity}
 \Hom^r_m\longrightarrow_{\BT}\Hom^r_n
\end{equation}
exists.  Apply Lemma~\ref{lem:all-arity-normalization}, and choose a
real \(W\) coding Turing functionals for the two continuous maps.  We
write \(\le_T\) for Turing reducibility and \(\oplus\) for the
computability-theoretic join of reals.  For every \(n\)-coloring \(c\),
its \(m\)-color query \(p\) satisfies
\[
 p\le_T W\oplus c.
\]
If \(A\) is an infinite exactly \(p\)-homogeneous set, then
\(B=\gamma(A)\le_TW\oplus A\) is \(c\)-almost homogeneous.  For some
finite \(N\), the set \(B\cap[N,\infty)\) is exactly homogeneous and
\begin{equation}\label{eq:cutoff-computation}
 B\cap[N,\infty)\le_TW\oplus A.
\end{equation}
There is no need to compute \(N\) uniformly: for each solution \(A\),
the relevant finite value can be built into the Turing program.  Thus
\eqref{eq:forbidden-fixed-arity} would yield, relative to \(W\),
\begin{equation}\label{eq:ramsey-computable-reduction}
 \mathrm{RT}^r_n\le_{\mathrm{sc}}^W\mathrm{RT}^r_m,
 \quad\text{and hence}\quad
 \mathrm{RT}^r_n\le_{\mathrm c}^W\mathrm{RT}^r_m.
\end{equation}
Here \(\le_{\mathrm{sc}}^W\) denotes strong computable reducibility
relative to \(W\), whereas \(\le_{\mathrm c}^W\) denotes computable
reducibility relative to \(W\). The latter follows from the former
because an ordinary computable reduction may also use the original
instance when computing a solution.
Patey's color-separation theorem states that
\[
 \mathrm{SRT}^r_n\not\le_{\mathrm c}\mathrm{RT}^r_m
 \qquad(r\ge2,\ n>m\ge2),
\]
and its proof relativizes to every oracle.  This contradicts
\eqref{eq:ramsey-computable-reduction}; see
\cite[Theorem~3.13 and Corollary~3.14]{Patey}.  The use of computable,
rather than merely Weihrauch, nonreducibility is what permits the
nonuniform removal of the finite exceptional set.

For \(r=1\), a direct selector argument gives the same separation.
Let the target coloring be \(c(k)=k\bmod n\), and let
\(p:\N\to\NatLt{m}\) be any query.  In each \(c\)-color class choose an
infinite subset on which \(p\) is constant.  Since \(n>m\), two such
sets \(U,V\), belonging to distinct \(c\)-colors, have the same
\(p\)-color.  Build an input by alternately appending a finite block
from \(U\) and from \(V\), extending each block until the online
selector chooses a new point.  Each stage terminates because the
selector preserves infinitude, and online decisions persist under
later extensions.  The resulting input is \(p\)-homogeneous, whereas
its selected output contains infinitely many points of each of two
different \(c\)-colors.  It is therefore not \(c\)-almost homogeneous.

We have proved, for every \(r\ge1\),
\begin{equation}\label{eq:fixed-arity-separation}
 n>m\ge2
 \quad\Longrightarrow\quad
 \Hom^r_m\not\longrightarrow_{\BT}\Hom^r_n.
\end{equation}

\subsection{Arity cannot be increased}

Finally suppose \(r<s\) and that
\(\Hom^r_m\to_{\BT}\Hom^s_n\).  The positive reduction already
proved, applied with source arity \(s\) and target arity \(r\), gives
\[
 \Hom^s_n\longrightarrow_{\BT}\Hom^r_{m+1}.
\]
Composing the two morphisms would give
\[
 \Hom^r_m\longrightarrow_{\BT}\Hom^r_{m+1},
\]
contrary to \eqref{eq:fixed-arity-separation}.  Together with the
positive reductions, this proves Theorem~\ref{thm:complete-classification}.

\begin{remark}\label{rem:one-color-boundary}
If one-color relations are included, the boundary cases are immediate.
For \(n=1\) a morphism always exists. If \(m=1<n\), suppose toward a
contradiction that \((\alpha,\beta)\) is a morphism. Fix \(A\in\X\) and
write \(\beta(A)=\{b_0<b_1<\cdots\}\). Define a target coloring
\(c\in\Ccal^s_n\) on \([\beta(A)]^s\) by
\[
 c(t)=i\bmod 2\qquad\text{when }\min t=b_i,
\]
and extend it arbitrarily to \([\N]^s\). The unique one-color source
coloring \(\alpha(c)\) is homogeneous on \(A\), but every tail of
\(\beta(A)\) contains \(s\)-sets of both \(c\)-colors. Hence
\(\beta(A)\) is not almost homogeneous for \(c\), contradicting
\eqref{eq:one-stage}. Thus
Theorem~\ref{thm:complete-classification} contains all nontrivial
finite-color cases.
\end{remark}

\begin{remark}\label{rem:complete-classification-scope}
The normalization proof uses Borelness of the challenge map when it
passes to a comeager continuity set. Thus the theorem is a
Borel--Tukey classification; it makes no assertion about arbitrary
challenge maps.
\end{remark}

\section{\texorpdfstring{
\(\Hom^2_2;\Hom^2_2;\Hom^2_2
 \longrightarrow_{\BT}\Hom^3_2\)}
{Three binary stages for binary triples}}
\label{sec:three-stage}

We construct the three-stage morphism by stabilizing sections,
homogenizing their limits, and then controlling stabilization rates.

\begin{lemma}[Lexicographically monotone sequences]\label{lem:lexconvergence}
Let \(k\in\N\) with \(k\ge1\). Every lexicographically nondecreasing or nonincreasing sequence in \((\NatLt{k})^\N\) converges in the product topology. The same conclusion holds if the sequence is only eventually monotone.
\end{lemma}
\begin{proof}
For a nondecreasing sequence \((u_j)_j\), the zeroth coordinates are
eventually constant. Once the coordinates below \(t\) are constant,
the \(t\)-th coordinates are nondecreasing on that tail and hence also
eventually constant. Induction gives convergence in \((\NatLt{k})^\N\). The
nonincreasing case is symmetric, and eventual monotonicity follows by
discarding an initial segment.
\end{proof}

\begin{theorem}\label{thm:three-stage}
There is a Borel morphism
\[
 \Hom^2_2;\Hom^2_2;\Hom^2_2
        \longrightarrow_{\BT}\Hom^3_2.
\]
\end{theorem}
\begin{remark}[Computability-theoretic comparison]
This Borel upper bound is compatible with a finite-iteration separation
of Patey and Shafer~\cite[Theorem~1.7]{PateyShafer}.  They prove that,
for every \(k\in\N\),
\[
 \mathsf{BRT}^{3}_{2,4}
 \not\leq_{\N}^{k}\mathsf{RT}^{2}_{2},
\]
where \(\mathsf{BRT}^{3}_{2,4}\) restricts the binary triple colorings
to those with no color-\(1\) homogeneous set of size \(4\).  Since an
infinite homogeneous set for such a coloring must have color \(0\),
their theorem in particular rules out solving
\(\mathsf{RT}^{3}_{2}\) by any fixed finite number of applications of
\(\mathsf{RT}^{2}_{2}\) in this computability-theoretic sense.  In the
reduction game defining \(\leq_{\N}^{k}\), successive instances and
the final solution must be Turing-computable from the initial instance
and the preceding solutions.  The stage maps in
Theorem~\ref{thm:three-stage} are only required to be Borel; the limit
and stabilization tests in \eqref{eq:query-two-total} and
\eqref{eq:stabilization-bound} need not be computable.  Thus the two
results concern different regularity requirements and do not conflict.
\end{remark}
\begin{proof}
Fix an enumeration \([\N]^2=\{s_t:t\in\N\}\). For a target
\(c\in\Ccal^3_2\), define the three queries as follows. Each formula
is total, including on invalid responses.

\emph{First query: stabilize all sections.}
For \(z\in\N\), define \(u^c_z\in\Bits^\N\) by
\[
 u^c_z(t)=\begin{cases}
 c(s_t\cup\{z\}),&z>\max s_t,\\
 0,&z\le\max s_t.
 \end{cases}
\]
Set
\begin{equation}\label{eq:query-one}
 \alpha_1(c)(x<y)=0
 \quad\Longleftrightarrow\quad u^c_x\le_{\rm lex}u^c_y.
\end{equation}
Suppose \(A=\{a_0<a_1<\cdots\}\) is almost homogeneous for this query. If its eventual color is zero, \((u^c_{a_j})_j\) is eventually nondecreasing. If that color is one, it is eventually strictly decreasing. By Lemma~\ref{lem:lexconvergence}, every coordinate is eventually constant. Therefore
\begin{equation}\label{eq:sectionlimits}
 \lim_{q\to\infty}c(x,y,a_q)
 \quad\text{exists for every fixed }x<y.
\end{equation}
Only sufficiently large \(q\), with \(a_q>y\), are relevant in this notation.

\emph{Second query: homogenize the limit pair coloring.}
For every \(c,A\), define a binary pair coloring \(\lambda_{c,A}\) by the total formula
\begin{equation}\label{eq:query-two-total}
 \lambda_{c,A}(i,j)=1
 \quad\Longleftrightarrow\quad
 (\exists N>j)(\forall q\ge N)\ c(a_i,a_j,a_q)=1
 \qquad(i<j).
\end{equation}
When the sequence fails to converge, this formula assigns the default value zero. For a valid first response, \eqref{eq:sectionlimits} gives
\[
 \lambda_{c,A}(i,j)=\lim_q c(a_i,a_j,a_q).
\]
Take \(\alpha_2(c,A)=\lambda_{c,A}\), and let
\(B=\{b_0<b_1<\cdots\}\) be an almost homogeneous response. The
limit colors are now constant on a tail of \(B\), but their
stabilization times may vary. The third query controls them.

\emph{Third query: exceed the stabilization bounds.}
For \(i<j\), if the following set is nonempty, let \(\tau_{c,A,B}(i,j)\) be its least element:
\begin{equation}\label{eq:stabilization-bound}
 \left\{N>b_j:(\forall q\ge N)\
 c(a_{b_i},a_{b_j},a_q)=\lambda_{c,A}(b_i,b_j)\right\}.
\end{equation}
If it is empty, set \(\tau_{c,A,B}(i,j)=b_j+1\). For valid \(A\), the set is always nonempty, for every \(B\). Put
\begin{equation}\label{eq:maximum-bound}
 T_{c,A,B}(j)=\max\bigl(\{0\}\cup
             \{\tau_{c,A,B}(i,j):i<j\}\bigr).
\end{equation}
The third query is
\begin{equation}\label{eq:query-three}
 \alpha_3(c,A,B)(j<k)=0
 \quad\Longleftrightarrow\quad b_k\ge T_{c,A,B}(j).
\end{equation}
This query has no infinite color-\(1\) homogeneous set: after fixing
its first member \(j\), every sufficiently large \(k\) satisfies
\(b_k\ge T_{c,A,B}(j)\). Hence every valid response \(C\) is
eventually homogeneous in color zero.

\emph{Verification of the output.}
Use the response
\begin{equation}\label{eq:final-response}
 \Gamma(A,B,C)=e_A[e_B[C]].
\end{equation}
Assume all three responses are valid. Choose a cutoff beyond which \(B\) is homogeneous for \(\lambda_{c,A}\), with color \(\varepsilon\). Since \(b_j\to\infty\), sufficiently large indices \(i<j\) have \(\lambda_{c,A}(b_i,b_j)=\varepsilon\). Also take \(i<j<k\) sufficiently late in \(C\) that the third query is zero on their pairs. Then
\[
 b_k\ge T_{c,A,B}(j)\ge\tau_{c,A,B}(i,j),
\]
so
\begin{equation}\label{eq:output-homogeneity}
 c(a_{b_i},a_{b_j},a_{b_k})
   =\lambda_{c,A}(b_i,b_j)=\varepsilon.
\end{equation}
Since \(j\mapsto a_{b_j}\) is strictly increasing, sufficiently late
triples in \(C\) form a tail of \(e_A[e_B[C]]\). This proves
\eqref{eq:three-stage-condition}, and \eqref{eq:final-response} uses
only the responses.

\emph{Borel regularity and totality.}
Evaluation at finitely many entries of increasing enumerations is
continuous. The lexicographic order on \(\Bits^\N\) is closed, so
\(\alpha_1\) is Borel. Formula \eqref{eq:query-two-total} makes each
coordinate of \(\alpha_2\) Borel. In
\eqref{eq:stabilization-bound}, the condition for a fixed \(N\) is
Borel; so are the least-witness operation, finite maxima, and the
comparison in \eqref{eq:query-three}. Thus all three stage maps are
jointly Borel on their full domains.

The image map \((A,D)\mapsto e_A[D]\) is continuous: membership of
\(n\) depends only on \(n\in A\) and \(\abs{A\cap\NatLt{n}}\in D\).
Applying it twice proves continuity of \(\Gamma\).
\end{proof}

\begin{remark}[Why a limit coloring alone is insufficient]\label{rem:rates}
The stable triple coloring
\[
 c(x<y<z)=1\quad\Longleftrightarrow\quad z=y+1
\]
has \(\lim_z c(x,y,z)=0\) for every \(x<y\), but every tail of
\(\N\) has triples of both colors. Sectional limits and a
homogeneous limit coloring therefore do not suffice. The third query
controls their stabilization times uniformly.
\end{remark}

\section{Two-stage reductions}\label{sec:two-stage}

We give two positive constructions: a three-color second stage for
binary triples, and two binary stages for arbitrary finite pair
palettes.

\subsection{\texorpdfstring{
\(\Hom^2_2;\Hom^2_3\longrightarrow_{\BT}\Hom^3_2\)}
{A binary stage followed by a three-color stage}}

The last two queries of Theorem~\ref{thm:three-stage} can be combined
into a single pair query with three colors.

\begin{theorem}\label{thm:two-stage-triples}
There is a Borel morphism
\[
 \Hom^2_2;\Hom^2_3\longrightarrow_{\BT}\Hom^3_2.
\]
\end{theorem}
\begin{proof}
Use the response \(\Gamma(A,B)=e_A[B]\).
Given \(c\in\Ccal^3_2\), use the first query
\(\alpha_1(c)\) from \eqref{eq:query-one}. For every \(A\in\X\),
write \(a_i=e_A(i)\), and let \(\lambda_{c,A}\) be the total Borel
pair coloring from \eqref{eq:query-two-total}. If \(A\) is a valid
first response, then
\[
 \lambda_{c,A}(i,j)=\lim_{q\to\infty}c(a_i,a_j,a_q)
 \qquad(i<j).
\]

For all \(c,A\) and \(i<j\), define \(\sigma_{c,A}(i,j)\) to be the
least member of
\begin{equation}\label{eq:two-stage-rate}
 \{N>j:(\forall q\ge N)\
       c(a_i,a_j,a_q)=\lambda_{c,A}(i,j)\},
\end{equation}
if this set is nonempty, and put \(\sigma_{c,A}(i,j)=j+1\)
otherwise. Set
\begin{equation}\label{eq:two-stage-gap}
 F_{c,A}(j)=
 \max\bigl(\{j+1\}\cup\{\sigma_{c,A}(i,j):i<j\}\bigr).
\end{equation}
These definitions remain total when the first response is invalid.
The second query is the three-color pair coloring
\begin{equation}\label{eq:two-stage-three-color}
 q_{c,A}(j,k)=
 \begin{cases}
  \lambda_{c,A}(j,k),&k\ge F_{c,A}(j),\\
  2,&k<F_{c,A}(j)
 \end{cases}
 \qquad(j<k).
\end{equation}
Take \(\alpha_2(c,A)=q_{c,A}\).

There is no infinite homogeneous set of color \(2\) for \(q_{c,A}\):
after fixing its first point \(j\), any sufficiently large later point
\(k\) satisfies \(k\ge F_{c,A}(j)\). Thus, if
\(q_{c,A}H^2_3B\), some tail \(B\cap[N,\infty)\) is homogeneous in a
color \(\varepsilon\in\Bits\).

Assume also that \(A\) is a valid first response. Then all the sets
in \eqref{eq:two-stage-rate} are nonempty. For
\(i<j<k\) in \(B\cap[N,\infty)\), applying
\eqref{eq:two-stage-three-color} to \((i,j)\) and \((j,k)\) gives
\[
 \lambda_{c,A}(i,j)=\varepsilon,\qquad
 k\ge F_{c,A}(j)\ge\sigma_{c,A}(i,j).
\]
Consequently,
\[
 c(a_i,a_j,a_k)=\lambda_{c,A}(i,j)=\varepsilon.
\]
Since \(e_A\) is increasing, this proves
\(cH^3_2e_A[B]\).

For each \(N\), the condition in
\eqref{eq:two-stage-rate} is Borel in \((c,A)\).
Least witnesses with the stated default, finite maxima, and the
comparison in \eqref{eq:two-stage-three-color} are Borel operations.
Thus the second query is jointly Borel on its full domain.
The first query is Borel by Theorem~\ref{thm:three-stage}, and
\((A,B)\mapsto e_A[B]\) is continuous.
\end{proof}

\begin{remark}\label{rem:exceptional-color}
In \eqref{eq:two-stage-three-color}, color \(2\) only marks failure of
the gap condition. Its graph is locally finite: each vertex has only
finitely many color-\(2\) neighbors. Thus the construction uses a
restricted kind of three-color query. One-stage palette separation
alone does not show that this third color is needed here.
\end{remark}

\subsection{\texorpdfstring{
\(\Hom^2_2;\Hom^2_2\longrightarrow_{\BT}\Hom^2_k\)}
{Two binary stages for finite pair palettes}}

\begin{theorem}\label{thm:two-stage-palettes}
For every finite positive integer \(k\), there is a Borel morphism
\[
 \Hom^2_2;\Hom^2_2\longrightarrow_{\BT}\Hom^2_k.
\]
\end{theorem}
\begin{proof}
Use the response \(\Gamma(A,B)=e_A[B]\).
Fix \(d\in\Ccal^2_k\). For \(z\in\N\), define \(v^d_z\in (\NatLt{k})^\N\)
by
\[
 v^d_z(t)=
 \begin{cases}
  d(t,z),&t<z,\\
  0,&t\ge z.
 \end{cases}
\]
The first binary query is
\begin{equation}\label{eq:pair-palette-first}
 p_1(d)(x,y)=0
 \quad\Longleftrightarrow\quad
 v^d_x\le_{\mathrm{lex}}v^d_y
 \qquad(x<y).
\end{equation}
If \(p_1(d)H^2_2A\), where \(a_i=e_A(i)\), the sequence
\((v^d_{a_i})_i\) is eventually lexicographically nondecreasing or
strictly decreasing. Lemma~\ref{lem:lexconvergence} therefore gives
the existence of \(\lim_q d(x,a_q)\) for every fixed \(x\), with only
the terms \(a_q>x\) considered.

For every \(d,A\) and \(i\), define \(\ell_{d,A}(i)\) to be the
unique \(\varepsilon\in\NatLt{k}\) satisfying
\begin{equation}\label{eq:pair-palette-limit}
 (\exists N>i)(\forall q\ge N)\
 d(a_i,a_q)=\varepsilon,
\end{equation}
if such an \(\varepsilon\) exists; otherwise set
\(\ell_{d,A}(i)=0\).
Define \(\rho_{d,A}(i)\) to be the least \(N>i\) such that
\[
 (\forall q\ge N)\ d(a_i,a_q)=\ell_{d,A}(i),
\]
if one exists, and set \(\rho_{d,A}(i)=i+1\) otherwise.
For a valid first response, these are the actual row limits and
their stabilization bounds.

The second binary query is
\begin{equation}\label{eq:pair-palette-second}
 p_2(d,A)(i,j)=0
 \quad\Longleftrightarrow\quad
 \ell_{d,A}(i)=\ell_{d,A}(j)
 \ \land\ j\ge\rho_{d,A}(i)
 \qquad(i<j).
\end{equation}
This query has no infinite homogeneous set of color \(1\).
Indeed, if \(C\) were such a set, the finiteness of \(k\) would give
an infinite \(C'\subseteq C\) on which \(\ell_{d,A}\) is constant.
Fix \(i\in C'\) and choose \(j\in C'\) with
\(j>i\) and \(j\ge\rho_{d,A}(i)\).
Then \eqref{eq:pair-palette-second} gives color \(0\), a contradiction.

Now suppose both responses \(A,B\) are valid.
Choose \(N\) such that \(p_2(d,A)\) is constantly \(0\) on
\([B\cap[N,\infty)]^2\). Equation~\eqref{eq:pair-palette-second}
implies that \(\ell_{d,A}\) is constant, with some value
\(\varepsilon\in\NatLt{k}\), on \(B\cap[N,\infty)\).
For \(i<j\) in this tail, it also gives \(j\ge\rho_{d,A}(i)\).
Hence
\[
 d(a_i,a_j)=\ell_{d,A}(i)=\varepsilon.
\]
Thus \(e_A[B]\) is almost homogeneous for \(d\).

The lexicographic order on \((\NatLt{k})^\N\) is closed, so the first
query is Borel. Formula~\eqref{eq:pair-palette-limit} gives Borel
coordinate maps for \(\ell_{d,A}\); taking the least stabilization
bound with the stated default gives Borel coordinate maps for
\(\rho_{d,A}\). Thus the second query is jointly Borel and total.
The response map is continuous and uses only \(A,B\).
\end{proof}


\section{\texorpdfstring{
\(\Hom^2_k;\Hom^2_2\not\longrightarrow_{\BT}\Hom^3_2\)}
{No finite-color first pair stage followed by a binary pair stage for binary triples}}
\label{sec:two-normalization}\label{sec:rapid-impossibility}

Throughout this section, fix a finite positive integer \(k\).
We prove that no such two-stage reduction exists. We first put the first
query in rapid form and the response in admissible form, then rule out
every reduction in this form; the second query need not be regular.


\subsection{Reduction setup}

A proposed two-stage reduction is given by
\begin{align*}
 \alpha_1&:\Ccal^3_2\to\Ccal^2_k
     &&\text{(the first query)},\\
 \alpha_2&:\Ccal^3_2\times\X\to\Ccal^2_2
     &&\text{(the second query)},\\
 \Gamma&:\X^2\to\X
     &&\text{(the final response map)}.
\end{align*}
Given a target coloring \(c\), a valid response \(A\) to
\(\alpha_1(c)\) determines the second query \(\alpha_2(c,A)\); a valid
response \(B\) to that query must make \(\Gamma(A,B)\) a valid response
to \(c\). Thus correctness means
\[
 \alpha_1(c)H^2_kA\ \land\ \alpha_2(c,A)H^2_2B
 \quad\Longrightarrow\quad cH^3_2\Gamma(A,B).
\]
We assume that \(\alpha_1\) and \(\Gamma\) are Borel, while
\(\alpha_2\) may be arbitrary.

For the transformed reduction produced by the normalization, we write
\[
 p:\Ccal^3_2\to\Ccal^2_k,\qquad
 q:\Ccal^3_2\times\X\to\Ccal^2_2,\qquad
 \theta:\Bits^\N\times\Bits^\N\to\Bits^\N.
\]
Here \(p,q\) are the normalized first and second queries, and
\(\theta\restriction\X^2\) is the normalized final response map. Their
correctness condition is the displayed implication above with
\((\alpha_1,\alpha_2,\Gamma)\) replaced by \((p,q,\theta)\). The first
half of the section makes \(p\) rapid and \(\theta\) admissible. The
second half proves the impossibility result without imposing any
regularity assumption on \(q\).

\begin{definition}\label{def:admissible}
A map \(\theta:\Bits^\N\times\Bits^\N\to\Bits^\N\) is
\emph{admissible} if, for all \(A,B\subseteq\N\) and \(n\in\N\),
\begin{align}
 n\in\theta(A,B)&\ \Longleftrightarrow\
 n\in\theta(A\cap\NatLt{n+1},B\cap\NatLt{n+1}),\label{eq:bi-prefix}\\
 \theta(A,B)&\subseteq A\cup B,\label{eq:bi-support}\\
 |\theta(A,B)\cap\NatLt{n+1}|&\le
 \min\{|A\cap\NatLt{n+1}|,|B\cap\NatLt{n+1}|\},\label{eq:bi-count}
\end{align}
and \(\theta(A,B)\) is infinite whenever both inputs are infinite.
We call \eqref{eq:bi-prefix} the prefix condition,
\eqref{eq:bi-support} the support condition, and
\eqref{eq:bi-count} the counting condition. Together, these
conditions make \(\theta\) continuous on the compact product.

A first query \(p:\Ccal^3_2\to\Ccal^2_k\) is \emph{rapid} if
\begin{equation}\label{eq:rapid-query}
 p(c)(i,j)=P_{i,j}(c\restr[\NatLt{j+1}]^3)\qquad(i<j)
\end{equation}
for some functions \(P_{i,j}:\Bits^{[\NatLt{j+1}]^3}\to\NatLt{k}\).
Thus a rapid query reads its value at \((i,j)\) from the target
coloring by stage \(j\), with no look-ahead.
\end{definition}

\subsection{Admissible responses $\theta$ instead of $\Gamma$}

We replace a Borel two-variable response $\Gamma$ by a continuous subresponse
with controlled behavior on finite inputs.

The space \(\mathcal W\) below compactly codes stagewise selection procedures. At
stage \(n\), a table reads only \(C\cap\NatLt{n+1}\) and either selects one
output point at or above \(n\), or records no selection by \(\infty\).

Let \(\Ninf=\N\cup\{\infty\}\) have its one-point
compactification topology. Consider the compact zero-dimensional space
\begin{equation}\label{eq:table-space}
 \begin{multlined}
 \mathcal W=\bigl\{w\in\Ninf^{\{(n,t):n\in\N,\ t\subseteq\NatLt{n+1}\}}:\\
 \forall n\in\N\ \forall t\subseteq\NatLt{n+1}\quad
 \bigl[w(n,t)\ge n\ \land\ (n\notin t\Rightarrow w(n,t)=\infty)\bigr]\bigr\}.
 \end{multlined}
\end{equation}
Here \(\infty\ge n\) for every \(n\). For \(w\in\mathcal W\) and
\(C\subseteq\N\), put
\begin{equation}\label{eq:table-evaluation}
 S_w(C)=\{w(n,C\cap\NatLt{n+1}):n\in\N\}\setminus\{\infty\}.
\end{equation}
The map \((w,C)\mapsto S_w(C)\) is continuous into \(\Bits^\N\).
Indeed, \(v\in S_w(C)\) iff
\(w(n,C\cap\NatLt{n+1})=v\) for some \(n\le v\), a finite union of clopen
conditions since \(\{v\}\) is clopen in \(\Ninf\). For finite \(t\),
we also have
\[
 |S_w(t)|\le |t|.
\]

\begin{claim}\label{clm:one-variable-package}
For every Borel \(\beta:\X\to\X\), there are a strictly increasing
\(u:\N\to\N\) and a table \(w\in\mathcal W\) such that
\begin{equation}\label{eq:one-variable-package}
 S_w(C)\in[\beta(u[C])]^\N\qquad(C\in\X).
\end{equation}
\end{claim}
\begin{claimproof}
Apply Lemma~\ref{lem:tail-normal} to \(\beta\), and take
\(u,F,(Q_n)_n\) as there. Define \(w(n,t)\) as the least point of
\(Q_n\cap F(t)\), and put \(w(n,t)=\infty\) when the intersection is
empty. Since \(Q_n\subseteq[n,\infty)\), the sets \(Q_n\) are pairwise
disjoint, and
\(F(t)\subseteq\bigcup_{m\in t}Q_m\), we have \(w\in\mathcal W\).
For \(C\in\X\), the last identity from
Lemma~\ref{lem:tail-normal} gives
\[
 Q_n\cap F(C)=Q_n\cap F(C\cap\NatLt{n+1}).
\]
Consequently,
\[
 w(n,C\cap\NatLt{n+1})=
 \begin{cases}
  \min(Q_n\cap F(C)),&Q_n\cap F(C)\ne\varnothing,\\
  \infty,&Q_n\cap F(C)=\varnothing.
 \end{cases}
\]
Thus \(S_w(C)\) contains exactly the least point of each \(Q_n\) that
meets \(F(C)\). Since \(F(C)\) is infinite and the sets \(Q_n\) are
finite, infinitely many \(Q_n\) meet \(F(C)\). Therefore
\(S_w(C)\in[F(C)]^\N\subseteq[\beta(u[C])]^\N\).
\end{claimproof}

We now isolate the local uniformization used below. Its proof uses
continuous reading of Mathias names, Silver's Ramsey theorem, and
Shoenfield absoluteness. See \cite[Sections~1 and~4]{Sabok}; the
continuous-reading observation follows Lemma~7.

\begin{lemma}[Coanalytic uniformization on a Mathias condition]
\label{lem:mathias-uniformization}
Let \(Z\) be a Polish space and let \(R\subseteq\X\times Z\) be
coanalytic. The definition of \(R\) may use fixed real parameters.
Suppose that, from stated properties of these parameters that are
preserved by forcing, ZFC proves uniformly that every vertical section
of \(R\) is nonempty:
\[
 (\forall a\in\X)(\exists z\in Z)\ R(a,z).
\]
For example, a parameter may be a code for a Borel map.
Then there are a set
\[
 [s,H]=\{s\cup A:A\in[H]^\N\},\qquad s<H,
\]
and a continuous map \(f:[s,H]\to Z\) such that
\[
 R(a,f(a))\qquad(a\in[s,H]).
\]
\end{lemma}
\begin{proof}
Let \(\dot a\) be the canonical name for the Mathias generic real.
The same ZFC proof of totality is valid in every Mathias forcing
extension, and hence
\[
 \Vdash(\exists z\in Z)\ R(\dot a,z).
\]
Choose a name \(\dot z\) such that
\[
 \Vdash R(\dot a,\dot z).
\]
By continuous reading of Mathias names, after strengthening a
condition there are \([s,H]\) and a continuous map
\(f:[s,H]\to Z\) such that
\[
 (s, H) \Vdash R(\dot a,f(\dot a)).
\]
The set
\[
 E=\{a\in[s,H]:\neg R(a,f(a))\}
\]
is analytic. By Silver's theorem there is an infinite
\(H'\subseteq H\) such that \([s,H']\subseteq E\) or
\([s,H']\cap E=\varnothing\). The first alternative is impossible:
the assertion \([s,H']\subseteq E\) is \(\boldsymbol{\Pi}^1_2\), and hence remains
true by Shoenfield absoluteness in a Mathias extension below
\([s,H']\). The generic real would then belong to \(E\), contrary to
the forced statement above. Therefore \([s,H']\cap E=\varnothing\),
and \(f\restr[s,H']\) is the required uniformization.
\end{proof}

\begin{lemma}[Compact two-variable subresponse]\label{lem:compact-two}
For every Borel \(\Gamma:\X^2\to\X\), there are
\(s\in[\N]^{<\infty}\), a strictly increasing map
\(v:\N\to\N\) with \(s<v[\N]\), a Borel map
\(\X\to\N^\N\), \(A\mapsto u_A\), whose values are strictly
increasing, and a continuous map \(K:\Bits^\N\times\Bits^\N\to\Bits^\N\)
such that
\begin{equation}\label{eq:compact-two-subresponse}
 K(A,B)\in[\Gamma(s\cup v[A],u_A[B])]^\N\qquad(A,B\in\X),
\end{equation}
and, for every finite \(t\) and every \(D\subseteq\N\),
\begin{equation}\label{eq:finite-axis-bounds}
 |K(t,D)|\le|t|,\qquad |K(D,t)|\le|t|.
\end{equation}
\end{lemma}
\begin{proof}
For \(a\in\X\), let \(R(a,u,w)\) say that \(u\) is strictly
increasing, \(w\in\mathcal W\), and
\[
 \forall C\in\X\quad S_w(C)\in[\Gamma(a,u[C])]^\N.
\]
Relative to a Borel code for \(\Gamma\), the relation is coanalytic,
and its sections are nonempty by \eqref{eq:one-variable-package}. The
proof of Claim~\ref{clm:one-variable-package} is a ZFC proof uniform
in this code. Thus the code
for \(\Gamma\) is an allowed parameter in
Lemma~\ref{lem:mathias-uniformization}.
Consequently, Lemma~\ref{lem:mathias-uniformization}, applied with
\(Z=\N^\N\times\mathcal W\), gives a Mathias condition \([s,H]\) and
a continuous map
\[
 f=(f_0,f_1):[s,H]\longrightarrow\N^\N\times\mathcal W
\]
such that \(R(a,f_0(a),f_1(a))\) for every \(a\in[s,H]\).

\smallskip
\noindent\emph{Extend continuously.}
Put \(u(a)=f_0(a)\) and \(w_a=f_1(a)\). Apply
Lemma~\ref{lem:compact}, after identifying \(H\) with \(\N\) and
embedding \(\mathcal W\) as a closed subspace of \(\Bits^\N\), to the
Borel map
\[
 [H]^\N\longrightarrow\mathcal W,
 \qquad A\longmapsto w_{s\cup A}.
\]
After shrinking \(H\) to a further infinite set \(T\subseteq H\), the
lemma gives a continuous extension on \(\Bits^T\). Put \(v=e_T\) and
reindex this extension to obtain a continuous map
\[
 \Bits^\N\longrightarrow\mathcal W,\qquad
 A\longmapsto\bar w_A,
\]
such that \(\bar w_A=w_{s\cup v[A]}\) for every infinite \(A\). Put
\[
 K_0(A,B)=S_{\bar w_A}(B),\qquad u_A=u(s\cup v[A])\quad(A\in\X).
\]
The first map is jointly continuous, is infinite on pairs of infinite
inputs, and has the required subresponse property. It already satisfies
\(|K_0(D,t)|\le|t|\) for finite second inputs. Finally set
\begin{equation}\label{eq:two-rank-trim}
 K(A,B)=\{z\in K_0(A,B):|K_0(A,B)\cap\NatLt{z}|\in A\}.
\end{equation}
The new map is continuous and remains infinite on pairs of infinite
inputs. To verify the first bound in \eqref{eq:finite-axis-bounds}, fix
a finite \(t\subseteq\N\) and \(D\subseteq\N\). For \(z\in K_0(t,D)\),
the rank map
\[
 z\longmapsto |K_0(t,D)\cap\NatLt{z}|
\]
is injective. Definition~\eqref{eq:two-rank-trim} retains \(z\) only
when this rank belongs to \(t\), so \(|K(t,D)|\le |t|\). For the
second bound, if \(t\) is finite and \(D\subseteq\N\), then
\(K(D,t)\subseteq K_0(D,t)\), while the preceding bound for the table
selector gives \(|K_0(D,t)|\le |t|\). Hence
\(|K(D,t)|\le |t|\), as required.
\end{proof}


\begin{lemma}\label{lem:admissible-normalization}
Suppose a reduction
\(
 \Hom^2_k;\Hom^2_2\longrightarrow_{\BT}\Hom^3_2
\)
has Borel first query and Borel response. It can be replaced by one
whose response is admissible and whose first query is Borel. If the
original second query is Borel, so is the new one; otherwise no
regularity is needed of either second query.
\end{lemma}
\begin{proof}
Take \(s,v,(u_A)_A,K\) from Lemma~\ref{lem:compact-two}.

\smallskip
\noindent\emph{Normalize the response.}
Recursively
choose \(x_0<x_1<\cdots\) with \(x_{n+1}>x_n+1\). Write
\(X_n=\{x_j:j\in\NatLt{n}\}\), and recursively define the finite sets
of output coordinates assigned at earlier stages by
\[
 D_0=\varnothing,\qquad D_{n+1}=D_n\cup Q_n,
 \quad\text{equivalently}\quad D_n=\bigcup_{m<n}Q_m.
\]
Put
\[
 E_n=D_n\cup\bigcup_{a,b\subseteq X_n}K(a,b).
\]
Continuity at the finitely many pairs \((a,b)\) permits choosing
\(x_n\) so large that
\begin{equation}\label{eq:two-fusion-preservation}
 K(a\cup C,b\cup D)\cap E_n=K(a,b)\cap E_n
\end{equation}
whenever \(a,b\subseteq X_n\) and \(C,D\subseteq[x_n,\infty)\).
Let \(Q_n\) collect the finite increments obtained by adjoining
\(x_n\) to the first input, the second input, or both:
\[
 Q_n=\bigcup_{a,b\subseteq X_n}
 \Bigl[
 \Bigl(
   K(a\cup\{x_n\},b)
   \cup K(a,b\cup\{x_n\})
   \cup K(a\cup\{x_n\},b\cup\{x_n\})
 \Bigr)\setminus K(a,b)
 \Bigr].
\]
Equation~\eqref{eq:two-fusion-preservation} shows that
\(Q_n\cap E_n=\varnothing\). Since
\[
 \bigcup_{m<n}Q_m=D_n\subseteq E_n,
\]
we have \(Q_n\cap Q_m=\varnothing\) whenever \(m<n\); hence the sets
\(Q_n\) are pairwise disjoint. Assign \(h(z)=x_n\) for \(z\in Q_n\).
Complete \(h\) by mapping all unassigned vertices injectively into
\(\N\setminus X\), where \(X=\{x_n:n\in\N\}\). Then \(h\) has finite
fibers and \(h^{-1}(\{x_n\})=Q_n\).

Put \(e=e_X\) and define
\begin{equation}\label{eq:admissible-response}
 \theta(A,B)=e^{-1}[h[K(e[A],e[B])]].
\end{equation}
Since \(K(a,b)\subseteq E_n\) for \(a,b\subseteq X_n\),
\eqref{eq:two-fusion-preservation} gives
\[
 K(a,b)\subseteq K(a\cup C,b\cup D)
 \qquad(C,D\subseteq[x_n,\infty)).
\]
Thus the outputs on successive finite prefixes along \(X\) form an
increasing sequence. Continuity of \(K\) then gives, for all \(A,B\),
\[
 K(e[A],e[B])=
 \bigcup_n K(e[A\cap\NatLt{n}],e[B\cap\NatLt{n}]).
\]

Since \(h^{-1}(\{x_n\})=Q_n\), we have
\[
 n\in\theta(A,B)
 \quad\Longleftrightarrow\quad
 K(e[A],e[B])\cap Q_n\ne\varnothing.
\]
Now \(Q_n\subseteq D_{n+1}\subseteq E_{n+1}\), so
\eqref{eq:two-fusion-preservation} at stage \(n+1\) yields
\begin{equation}\label{eq:qn-prefix-dependence}
 K(e[A],e[B])\cap Q_n
 =K(e[A\cap\NatLt{n+1}],e[B\cap\NatLt{n+1}])\cap Q_n.
\end{equation}
This proves \eqref{eq:bi-prefix}. If \(n\notin A\cup B\), the two
prefixes on the right of \eqref{eq:qn-prefix-dependence} are subsets
of \(X_n\). Their \(K\)-value is therefore contained in \(E_n\),
whereas \(Q_n\cap E_n=\varnothing\). Hence the right-hand side is
empty, so \(n\notin\theta(A,B)\). This proves \eqref{eq:bi-support}.
Moreover,
\[
 \theta(A,B)\cap\NatLt{n+1}
 =e^{-1}\!\left[h\left[K(e[A\cap\NatLt{n+1}],e[B\cap\NatLt{n+1}])\right]\right].
\]
The inequalities \eqref{eq:finite-axis-bounds} give \eqref{eq:bi-count}.
If \(A\) and \(B\) are infinite, then \(e[A]\) and \(e[B]\) are
infinite, so \(K(e[A],e[B])\) is infinite by Lemma~\ref{lem:compact-two}.
Moreover, \(K(\varnothing,\varnothing)=\varnothing\) by
\eqref{eq:finite-axis-bounds}, and each increment between successive
finite-prefix outputs is contained in the corresponding \(Q_n\).
The union identity above therefore implies
\(K(e[A],e[B])\subseteq\bigcup_n Q_n\), so its image under \(h\)
is contained in \(X\). This image is infinite: otherwise
\(K(e[A],e[B])\) would be contained in a finite union of finite
fibers of \(h\). Since \(e\colon\N\to X\) is a bijection,
\(\theta(A,B)=e^{-1}[h[K(e[A],e[B])]]\) is infinite as well.

\smallskip
\noindent\emph{Transport the queries.}

The assumption and the conclusion we want are
\[
\alpha_1(\hat{c})H^2_k\hat{A}\ \land\ \alpha_2(\hat{c},\hat{A})H^2_2\hat{B}
\quad\Longrightarrow\quad \hat{c}H^3_2\Gamma(\hat{A},\hat{B}) \quad\text{for every $\hat{c}, \hat{A}, \hat{B}$},
\]
and
\[
p(c)H^2_kA\ \land\ q(c,A)H^2_2B
\quad\Longrightarrow\quad cH^3_2\theta(A,B) \quad\text{for every $c, A, B$}.
\]

So fix \(c\in\Ccal^3_2\) and $A, B \in [\N]^\N$.

Put \(\pi=e^{-1}\circ h:h^{-1}[X]\to\N\). Define \(\widehat c\in\Ccal^3_2\), for
\(z_0<z_1<z_2\), by
\[
 \widehat c(\{z_0,z_1,z_2\})=
 \begin{cases}
  c(\{\pi(z_0),\pi(z_1),\pi(z_2)\}),
   &\begin{array}{l}
     \text{if }h(z_0),h(z_1),h(z_2)\in X\text{ and}\\[-2pt]
     h(z_0),h(z_1),h(z_2)\text{ are pairwise distinct},
    \end{array}\\[6pt]
  0,&\text{otherwise}.
 \end{cases}
\]
Given original queries \(\alpha_1,\alpha_2\), use
\begin{align*}
 p(c)(i,j)&=\alpha_1(\widehat c)(v(e(i)),v(e(j))),\\
 q(c,A)(i,j)&=\alpha_2(\widehat c,s\cup v[e[A]])
                  (u_{e[A]}(e(i)),u_{e[A]}(e(j))).
\end{align*}
Valid new responses give the original responses
\[
 \hat{A} = s\cup v[e[A]]
 \quad\text{and}\quad
 \hat{B} = u_{e[A]}[e[B]].
\]
Hence \(K(e[A],e[B])\) is almost homogeneous for \(\widehat c\).
So, \(\theta(A,B)\) is almost homogeneous for \(c\). The formulas have the claimed regularity.
\end{proof}

\subsection{Rapid first queries $p$}

\begin{proposition}[Rapid normal form]\label{prop:rapid-normal-form}
If there is a reduction
\(
\Hom^2_k;\Hom^2_2\longrightarrow_{\BT}\Hom^3_2
\) with Borel first
query and Borel response, there is one with rapid first
query and admissible response. The second query remains Borel if it
was originally Borel; it may otherwise be arbitrary.
\end{proposition}
\begin{proof}
Apply Lemma~\ref{lem:admissible-normalization}, and call the resulting
maps \(p,q,\theta\).

\smallskip
\noindent\emph{Make the first query continuous.}
Choose a comeager \(G_\delta\) set \(D\) on which \(p\) is continuous,
and apply Lemma~\ref{lem:free-cube} to binary triple colorings to
obtain an infinite set \(Y\subseteq\N\) and a continuous map
\(\Psi:\Ccal^3_2\to D\) satisfying \eqref{eq:free-cube}. With
\(e=e_Y\), set
\begin{equation}\label{eq:free-cube-pullback}
 \begin{aligned}
  p_1(c)(i,j)&=p(\Psi(c))(e(i),e(j)),\\
  q_1(c,A)(i,j)&=q(\Psi(c),e[A])(e(i),e(j)),\\
  \theta_1(A,B)&=e^{-1}[\theta(e[A],e[B])].
 \end{aligned}
\end{equation}
The support condition \eqref{eq:bi-support} for \(\theta\) gives
\[
 \theta(e[A],e[B])\subseteq e[A]\cup e[B]\subseteq Y,
\]
so the old output lies inside \(Y\).
Equation~\eqref{eq:free-cube} therefore proves correctness of this
new reduction. The four admissibility conditions are preserved by this
increasing reindexing, and \(p_1\) is continuous on all of the compact
space \(\Ccal^3_2\).

\smallskip
\noindent\emph{Make the first query rapid.}
For each \(u<v\), continuity of \(c\mapsto p_1(c)(u,v)\) into the finite
discrete space \(\NatLt{k}\), together with compactness of \(\Ccal^3_2\),
gives a finite set
\(D_{u,v}\subseteq[\N]^3\) such that
\[
 c\restr D_{u,v}=c'\restr D_{u,v}
 \quad\Longrightarrow\quad
 p_1(c)(u,v)=p_1(c')(u,v).
\]
Choose \(z_0<z_1<\cdots\) so that
\(z_{n+1}\) is above all vertices in the sets
\(D_{z_i,z_j}\), \(i<j\le n\). Let \(Z=\{z_n:n\in\N\}\), and
define \(L_Z:\Ccal^3_2\to\Ccal^3_2\) by
\begin{equation}\label{eq:LZ-definition}
 L_Z(c)(t)=
 \begin{cases}
  c(e_Z^{-1}[t]),&t\in[Z]^3,\\
  0,&t\notin[Z]^3.
 \end{cases}
\end{equation}
Following the pullback pattern in \eqref{eq:free-cube-pullback}, set
\begin{equation}\label{eq:rapid-pullback}
 \begin{aligned}
  p_2(c)(i,j)&=p_1(L_Z(c))(e_Z(i),e_Z(j)),\\
  q_2(c,A)(i,j)&=q_1(L_Z(c),e_Z[A])(e_Z(i),e_Z(j)),\\
  \theta_2(A,B)&=e_Z^{-1}[\theta_1(e_Z[A],e_Z[B])].
 \end{aligned}
\end{equation}
The support condition \eqref{eq:bi-support} gives
\[
 \theta_1(e_Z[A],e_Z[B])
 \subseteq e_Z[A]\cup e_Z[B]\subseteq Z.
\]
If \(A,B\) are valid responses to the new queries, then \(e_Z[A]\)
and \(e_Z[B]\) are valid responses to the corresponding old queries.
Thus \(\theta_1(e_Z[A],e_Z[B])\) is almost homogeneous for \(L_Z(c)\).
Since this output lies in \(Z\), \eqref{eq:LZ-definition} shows that
its pullback \(\theta_2(A,B)\) is almost homogeneous for \(c\).
Hence \eqref{eq:rapid-pullback} defines a valid reduction, and the
increasing reindexing preserves admissibility.

Fix \(i<j\). By the choice of \(D_{z_i,z_j}\), the value
\(p_2(c)(i,j)\) depends only on
\(L_Z(c)\restr D_{z_i,z_j}\). Every vertex occurring in a triple in
\(D_{z_i,z_j}\) is below \(z_{j+1}\). If such a triple \(t\) is not
in \([Z]^3\), then \(L_Z(c)(t)=0\), independently of \(c\). If
\(t\in[Z]^3\), then
\[
 t=e_Z[s]\quad\text{for some }s\in[\NatLt{j+1}]^3,
 \qquad L_Z(c)(t)=c(s).
\]
Consequently, \(p_2(c)(i,j)\) depends only on
\(c\restr[\NatLt{j+1}]^3\). Thus \(p_2\) is rapid, as required.
\end{proof}

\subsection{Impossibility in rapid normal form}

By Proposition~\ref{prop:rapid-normal-form}, it suffices to rule out
reductions with a rapid first query and an admissible response. We do
this by constructing a target coloring that is nondecreasing in its last
coordinate, while allowing the second query to remain arbitrary.
Throughout the proof of Theorem~\ref{thm:rapid-impossibility}, we use
the convention \(\max\varnothing=-1\). Admissibility refers to the full
finite-input extension in Definition~\ref{def:admissible}.

\begin{theorem}[Uniform impossibility for the second query]\label{thm:rapid-impossibility}
For every finite positive integer \(k\), every rapid first query
\(p:\Ccal^3_2\to\Ccal^2_k\), and every
admissible \(\theta\), there are \(c\in\Ccal^3_2\) and an infinite exactly
\(p(c)\)-homogeneous set \(A\) such that
\begin{equation}\label{eq:rapid-1}
 \forall d\in\Ccal^2_2\ \exists B\in[A]^\N\quad
 \left\{\begin{array}{l}
 B\text{ is exactly }d\text{-homogeneous},\\
 \neg\bigl(cH^3_2\theta(A,B)\bigr).
 \end{array}\right.
\end{equation}
In particular, no second query \(q:\Ccal^3_2\times\X\to\Ccal^2_2\),
even an arbitrary one, completes \((p,\theta)\) to a valid two-stage
reduction.
\end{theorem}

The proof begins below and continues through
Subsection~\ref{subsec:rapid-second}, where \eqref{eq:rapid-1} is
established.

The final assertion of the theorem follows by taking \(d=q(c,A)\).
Recall that \(\le_T\)
denotes Turing reducibility and \(\oplus\) the computability-theoretic
join of reals. We use Galvin--Prikry and the following relativized
cone-avoidance theorem for pairs with \(k\) colors:
\begin{equation}\label{eq:pair-cone-avoidance}
 \begin{aligned}
 &\forall C\in\Bits^\N\ \forall W\in\Bits^\N\ \forall d\in\Ccal^2_k,\\
 &\bigl(C\not\le_T W\ \land\ d\le_T W\bigr)\quad\Longrightarrow\\
 &\qquad \exists H\in\X\,
 \bigl(H\text{ is exactly }d\text{-homogeneous and }C\not\le_T W\oplus H\bigr).
 \end{aligned}
\end{equation}
This is due to Seetapun and Slaman; see \cite[Section~1.1]{CP}.

\subsection{Fixing a Turng oracle $Z$ and its jump $Z'$}\label{subsec:rapid-responses}
For the remainder of the proof, fix a rapid first query
\(p:\Ccal^3_2\to\Ccal^2_k\) and an admissible response
\(\theta:\Bits^\N\times\Bits^\N\to\Bits^\N\), as in
Theorem~\ref{thm:rapid-impossibility}.
We thin the first-input space so that all responses to finite second
inputs vary continuously and can be encoded by one real.

Although the reduction uses \(\theta\) only on pairs of infinite
sets, admissibility gives a map on all of
\(\Bits^\N\times\Bits^\N\). Thus \(\theta(D,b)\) is already defined
when \(D\in\X\) and \(b\in[\N]^{<\infty}\). The counting condition
\eqref{eq:bi-count} implies, for every \(n\),
\[
 |\theta(D,b)\cap\NatLt{n+1}|
 \le |b\cap\NatLt{n+1}|\le |b|.
\]
Letting \(n\) increase gives
\begin{equation}\label{eq:rapid-2}
 |\theta(D,b)|\le |b|.
\end{equation}
Thus \(\theta(D,b)\) is a finite output. This estimate controls
only its cardinality, not the locations of its elements; in particular,
it gives no uniform bound on \(\max\theta(D,b)\) from \(b\)
alone. The thinning below makes the entire finite output, including
these locations, vary continuously with \(D\).

Give \(\mathrm{Fin}=[\N]^{<\infty}\) the discrete topology. For each
\(D\), define its finite-input response map by
\begin{equation}\label{eq:rapid-3}
\mathcal F(D)=\bigl(\theta(D,b)\bigr)_{b\in\mathrm{Fin}}
\in\mathrm{Fin}^{\mathrm{Fin}}
\end{equation}
The map \(D\mapsto\mathcal F(D)\) is Borel. For fixed \(b,n\), the assertion
\(n\in\theta(D,b)\) depends on a finite prefix of \(D\);
equality with a specified finite output set is therefore Borel.

\begin{lemma}[Ramsey continuity]\label{lem:ramsey-continuity}
If \(F:[\N]^\N\to Y\) is Borel and \(Y\) has a countable
base, there is an infinite \(H\) such that
\(F\restriction[H]^\N\) is continuous.
\end{lemma}
\begin{proof}
This is analogous to the fusion argument in Lemma~\ref{lem:compact}.
Enumerate the inverse images of a countable base of
\(Y\) as \(D_0,D_1,\ldots\). Construct a fusion sequence. At stage
\(n\), after choosing finitely many points of \(H\), apply
Galvin--Prikry successively to each of the finitely many pairs
consisting of a subset \(s\) of those points and an index \(\ell\le n\).
Thin the candidate set so that the corresponding cone is either contained
in \(D_\ell\) or disjoint from it.

On the diagonal set \(H\), each \(D_\ell\cap[H]^\N\) is clopen. A
long enough finite prefix determines the decision made during the
fusion. Hence the restriction of \(F\) is continuous.
\end{proof}

Apply the lemma to \eqref{eq:rapid-3}. Fix a real \(Z\) coding:

\begin{enumerate}
\item \(H\) and codes witnessing that \(p\) is rapid and \(\theta\) is
  admissible;
\item a continuous code for
  \(\mathcal F\restriction[H]^\N\).
\end{enumerate}

More precisely, fix a canonical natural-number code \(\ulcorner
a\urcorner\) for each finite set \(a\). The continuous code contained
in \(Z\) yields fixed Turing functionals \(\Phi\) and \(\mu\) such
that, for every \(D\in[H]^\N\) and \(b\in\mathrm{Fin}\), both
computations below halt, with
\[
 \Phi^{Z\oplus D}(\ulcorner b\urcorner)
   =\ulcorner\theta(D,b)\urcorner,
 \qquad
 N:=\mu^{Z\oplus D}(\ulcorner b\urcorner).
\]
The number \(N\) so obtained satisfies
\[
 \forall E\in[H]^\N\quad
 E\cap\NatLt{N}=D\cap\NatLt{N}
 \quad\Longrightarrow\quad
 \theta(E,b)=\theta(D,b).
\]

We write \(Z'\) for the Turing jump of \(Z\).

\subsection{Defining the coloring $c$}\label{subsec:rapid-candidates}
We construct the monotone triple coloring \(c\) and an infinite
set \(X\), and then thin each row \(X\cap g^{-1}(i)\) to a homogeneous
candidate set \(H_i\). The first response \(A\) will be assembled from
these sets \(H_i\).

Let \((J_s)_{s\in\N}\) be a uniformly \(Z\)-computable sequence
of finite sets such that
\[
 J_s\subseteq J_{s+1}\quad\text{and}\quad Z'=\bigcup_{s\in\N}J_s.
\]
Put
\begin{equation}\label{eq:rapid-4}
 j(x,y,z)=1\iff (J_z\setminus J_y)\cap\NatLt{x+1}\ne\varnothing.
\end{equation}
We use four elementary facts.

\begin{itemize}
\item
  \(j\) is nondecreasing in its last coordinate.
\item
  The limiting pair coloring
  \[
   j_\infty(x,y)=\lim_{z\to\infty}j(x,y,z)
  \]
  is, for each fixed \(x\), eventually zero as \(y\to\infty\). In
  particular, \(j\) has no infinite homogeneous set of color 1.
\item
  If some \(e\le h\) enters \(J\) after stage \(M\), then
  \begin{equation}\label{eq:rapid-5}
   j_\infty(x,y)=1\qquad(h\le x<y\le M).
  \end{equation}
\item
  If \(Z'\not\le_T V\) and \(Z\le_TV\), then every total
  \(V\)-computable function \(f:\N\to\N\) satisfying \(f(h)\ge h\) for
  all \(h\) has the following property:
  \[
   \forall N\ \exists h\ge N\ \exists e\le h\quad
   e\text{ enters \(J\) after \(f(h)\)}.
  \]
  If this failed, there would be an \(N\) such that, for every \(h\ge N\),
  \[
   Z'\cap\NatLt{h+1}=J_{f(h)}\cap\NatLt{h+1}.
  \]
  Thus \(Z\oplus f\), together with the finite set
  \(Z'\cap\NatLt N\), would compute \(Z'\). Since
  \(Z\oplus f\le_TV\), this contradicts \(Z'\not\le_TV\).
\end{itemize}

For a labeling \(g:\N\to\N\), define
\begin{equation}\label{eq:rapid-6}
 c_g(x,y,z)=
 \begin{cases}
  1,& g(x)\ne g(y),\\
  j(x,y,z),& g(x)=g(y).
\end{cases}
\end{equation}
For each \(i\in\N\), call the fiber \(g^{-1}(i)\) the \(i\)-th row.
Thus \(c_g\) agrees with \(j\) when its first two points lie in the same
row, and it equals 1 when they lie in different rows. Write
\(r_g=p(c_g)\in\Ccal^2_k\).

For each \(i<t\), the defining property of a rapid query implies that
\(r_g(i,t)\) depends only on
\(c_g\restriction[\NatLt{t+1}]^3\). If
\((x,y,z)\in[\NatLt{t+1}]^3\), then \(x<y<z\le t\), so \(x,y<t\), and
\(c_g(x,y,z)\) depends on \(g\) only through \(g(x)\) and \(g(y)\).
Hence all the values \(r_g(i,t)\), \(i<t\), are determined by
\(g\restriction\NatLt t\) and do not depend on \(g(t)\). Consequently,
\(g(t)\) can be chosen after these values have been fixed.

\begin{claim}\label{clm:rapid-labeling}
We can choose a labeling \(g:\N\to\N\) so that
\begin{equation}\label{eq:rapid-7}
 Z'\not\le_T Z\oplus g,
\end{equation}
and so that, for every finite set \(F\subseteq\N\) and every map
\(\eta:F\to\NatLt{k}\), the candidate set
\begin{equation}\label{eq:rapid-8}
 R(F,\eta)=\{t\in H:t>\max F,\
 r_g(x,t)=\eta(x)\text{ for every }x\in F\}
\end{equation}
has the following property: if \(R(F,\eta)\) is infinite, then
\(R(F,\eta)\cap g^{-1}(i)\) is infinite for every \(i\in\N\).
\end{claim}

\begin{claimproof}
Work in the Baire space \(\N^\N\) of labelings, and write
\(R_g(F,\eta)\) for the set in \eqref{eq:rapid-8} obtained from
\(g\). For fixed \(F,\eta,i,N\), let
\[
 \mathcal B_{F,\eta,i,N}
 =\{g\in\N^\N:R_g(F,\eta)\text{ is infinite and }
       R_g(F,\eta)\cap g^{-1}(i)\cap[N,\infty)=\varnothing\}.
\]
We show that \(\mathcal B_{F,\eta,i,N}\) is nowhere dense. Let
\(\sigma\in\N^{<\N}\) and put \(M=\max(N,|\sigma|)\). If
\(R_h(F,\eta)\subseteq\NatLt M\) for every \(h\in[\sigma]\), then
\([\sigma]\) is disjoint from \(\mathcal B_{F,\eta,i,N}\). Otherwise,
choose \(h\in[\sigma]\) and \(t\in R_h(F,\eta)\) with \(t\ge M\).
Membership of \(t\) in \(R_g(F,\eta)\) depends only on
\(g\restriction t\), so
\(\tau=(h\restriction t)^\frown\langle i\rangle\) extends \(\sigma\) and
every \(g\in[\tau]\) has \(t\in R_g(F,\eta)\cap g^{-1}(i)\cap[N,\infty)\).
Hence \([\tau]\) is disjoint from \(\mathcal B_{F,\eta,i,N}\), proving
that this set is nowhere dense.

The standard relativized upper-cone theorem states that if
\(x\not\le_T y\), then
\(\{z\in 2^\N:x\le_T y\oplus z\}\) is meager
\cite[Lemma~4.1]{KumarShelah}; its category proof works unchanged in
Baire space. Hence \(\{g:Z'\le_T Z\oplus g\}\) is meager. The Baire
category theorem now gives a labeling \(g\) satisfying both
requirements.
\end{claimproof}

Fix such a labeling \(g\). Set \(c=c_g\in\Ccal^3_2\),
\(r=p(c)\in\Ccal^2_k\), and \(W_0=Z\oplus g\). An infinite set \(X\)
is \emph{\(r\)-cohesive} if, for every \(x\in\N\),
there is a color \(\delta_x\in\NatLt{k}\) such that
\(r(x,y)=\delta_x\) for all sufficiently large \(y\in X\).

\begin{claim}\label{clm:rapid-cohesive}
There is an infinite \(X\subseteq H\) such that:

\begin{itemize}
\item
  \(X\) is \(r\)-cohesive;
\item
  every \(X\cap g^{-1}(i)\) is infinite;
\item
  \(Z'\not\le_T W_0\oplus X\).
\end{itemize}
\end{claim}

\begin{claimproof}
We use a cone-avoiding form of the Rasiowa--Sikorski argument. For
finite \(F\subseteq\N\), \(\eta:F\to\NatLt{k}\), and \(m\in\N\), put
\[
 R_{F,\eta,m}=R(F,\eta)\cap[m,\infty).
\]
Let \(\mathbb P\) consist of the tuples \(P=(s;F,\eta,m)\), where
\(s\subseteq H\) is finite, \(\max s<m\), and the reservoir
\(R_P=R_{F,\eta,m}\) is infinite. The condition \(P\) represents
\[
 [P]=\{Y\in[H]^\N:s\subseteq Y\subseteq s\cup R_P\}.
\]
Order \(\mathbb P\) by reverse inclusion of these classes. In
particular, one may extend \(P\) by adding finitely many points of
\(R_P\) to its stem and then replacing \(R_P\) by an infinite
candidate-set tail contained in the old reservoir. The poset
\(\mathbb P\) is countable, and every reservoir is
\(W_0\)-computable.

We identify three families of dense subsets of \(\mathbb P\). First,
for \(x\in\N\), let
\[
 \mathcal C_x=\{(s;F,\eta,m)\in\mathbb P:x\in F\}.
\]
To see that \(\mathcal C_x\) is dense, take
\(P=(s;F,\eta,m)\notin\mathcal C_x\) and partition its reservoir into
\[
 R_P^\delta=\{y\in R_P:y>x\text{ and }r(x,y)=\delta\},
 \qquad \delta<k.
\]
Since the partition is finite, at least one part, say \(R_P^\delta\),
is infinite. After deleting a finite initial segment, it has the form
\[
 R(F\cup\{x\},\eta\cup\{(x,\delta)\})\cap[m',\infty),
\]
for some sufficiently large \(m'\), and therefore gives an extension
in \(\mathcal C_x\).

Every infinite reservoir has infinitely many points in each
\(g\)-row, by Claim~\ref{clm:rapid-labeling}. It follows that, for
\(i,n\in\N\),
\[
 \mathcal L_{i,n}
 =\{(s;F,\eta,m)\in\mathbb P:|s\cap g^{-1}(i)|\ge n\}
\]
is dense: add enough points of \(R_P\cap g^{-1}(i)\) to the stem and
then pass to a tail.

Finally, for each Turing functional \(\Phi_e\), let
\[
 \mathcal D_e
 =\{P\in\mathbb P:(\forall Y\in[P])
       \Phi_e^{W_0\oplus Y}\ne Z'\}.
\]
The usual computable-Mathias cone-avoidance argument shows that
\(\mathcal D_e\) is dense. Indeed, fix \(P=(s;F,\eta,m)\). Using
\(W_0\), define a partial function \(\psi\) by searching, on input
\(n\), for a finite \(u\subseteq R_P\) and a computation
\begin{equation}\label{eq:rapid-cone-search}
 \Phi_e^{W_0\oplus(s\cup u)}(n)\mathord\downarrow=b,
\end{equation}
and taking the first value \(b\) found. If \(\psi(n)\) is undefined for
some \(n\), then
\(\Phi_e^{W_0\oplus Y}(n)\) diverges for every \(Y\in[P]\), since any
convergent computation has finite use and would occur in the search.
Thus \(P\in\mathcal D_e\). Otherwise
\(\psi\le_TW_0\), so \eqref{eq:rapid-7} gives an \(n\) with
\(\psi(n)\ne Z'(n)\). Add to the stem a finite \(u\) witnessing
\eqref{eq:rapid-cone-search} and move the reservoir above the use
of the computation, obtaining an extension \(Q\). Then every
\(Y\in[Q]\) satisfies
\(\Phi_e^{W_0\oplus Y}(n)=\psi(n)\ne Z'(n)\), so
\(Q\in\mathcal D_e\). This proves the density of \(\mathcal D_e\).

Apply the Rasiowa--Sikorski lemma to the countable family
\[
 \{\mathcal C_x:x\in\N\}
 \cup\{\mathcal L_{i,n}:i,n\in\N\}
 \cup\{\mathcal D_e:e\in\N\}.
\]
Let \((P_\ell)_{\ell\in\N}\) be a descending sequence meeting every member
of this family, and let \(X\) be the union of its stems. Then
\(X\in[P_\ell]\) for every \(\ell\). Meeting every \(\mathcal L_{i,n}\)
makes \(X\) infinite and ensures that each
\(X\cap g^{-1}(i)\) is infinite. Meeting every \(\mathcal C_x\) makes
\(X\) \(r\)-cohesive, and meeting every
\(\mathcal D_e\) gives \(Z'\not\le_TW_0\oplus X\).
\end{claimproof}

Fix such a set \(X\), and put \(W=W_0\oplus X\). For each \(i\), apply
the cone-avoidance theorem \eqref{eq:pair-cone-avoidance} to \(r\) on
\(X\cap g^{-1}(i)\). This gives an infinite homogeneous set \(H_i\) of
color \(\varepsilon_i<k\) such that
\begin{equation}\label{eq:rapid-9}
 Z'\not\le_T W\oplus H_i.
\end{equation}
Because \(X\) is \(r\)-cohesive and \(H_i\subseteq X\), every \(x\in H_i\)
satisfies \(r(x,y)=\varepsilon_i\) for all sufficiently large \(y\in X\).
Since \(k\) is finite, choose \(\varepsilon<k\) such that
\(I=\{i\in\N:\varepsilon_i=\varepsilon\}\) is infinite.

Every finite \(r\)-\(\varepsilon\)-homogeneous stem from
\(\bigcup_{i\in I}H_i\) is therefore compatible with a cofinite tail
of each \(H_i\). We build the first answer $A$ from these tails.

We use only the individual statements \eqref{eq:rapid-9}; neither
their countable join nor the final first answer must avoid \(Z'\).

\subsection{Virtual norms $(\nu_i)_i$}\label{subsec:rapid-virtual}
In this subsection, for each \(i\in I\), we define a compact hereditary
family \(L_i\) of finite row blocks that, for suitable parameters
\(u\) and \(h\), do not force the appearance of a \(j\)-zero triple. We then
define an increasing closure hierarchy \((K_{i,n})_{n\in\N}\) and the
associated virtual norm \(\nu_i\) on finite subsets of \(H_i\); the norm
of a block is the first level of the hierarchy that contains it.
Union closure makes the bounded-norm sets into an ideal, while
Ramsey closure ensures that a high-norm block retains a high-norm
homogeneous subblock for any binary coloring
\(d\in\Ccal^2_2\) given later. We use
virtual first inputs in these definitions so that the norms are
available before the actual first response, and hence before the second
coloring $d$, is known.

Fix \(i\in I\). For \(h\in\N\) and \(u\subseteq H\cap\NatLt{h}\), define
the virtual first input
\begin{equation}\label{eq:rapid-10}
 V_{u,h}^i=u\cup(H_i\cap[h,\infty)).
\end{equation}
This belongs to \([H]^\N\). For \(F\subseteq H_i\cap[h,\infty)\),
let
\begin{equation}\label{eq:rapid-11}
 \begin{split}
 F\in L_{u,h}^i
 \quad\Longleftrightarrow\quad
 \forall a\in[F]^{<\infty}\ \exists s\subseteq u\quad&
 \theta(V_{u,h}^i,s\cup a)\cap[h,\infty)\\
 &\text{contains no \(j\)-zero triple}.
 \end{split}
\end{equation}
This family is hereditary and closed in the product topology on
\(2^{H_i\cap[h,\infty)}\).

\begin{claim}\label{clm:rapid-local-no-infinite}
\(L_{u,h}^i\) has no infinite member.
\end{claim}
\begin{claimproof}
Suppose that \(F\in L_{u,h}^i\) were infinite. For each
\(s\subseteq u\), both
\(V_{u,h}^i\) and \(s\cup F\) are infinite, so admissibility makes
\[
 R_s=\theta(V_{u,h}^i,s\cup F)
\]
infinite. The support condition \eqref{eq:bi-support} gives
\[
 R_s\subseteq V_{u,h}^i\cup(s\cup F)=V_{u,h}^i
   =u\cup(H_i\cap[h,\infty)).
\]
Hence \(R_s\cap[h,\infty)\subseteq H_i\). This tail is still infinite.
As \(j\) has no infinite homogeneous set of
color \(1\), choose a \(j\)-zero triple
\[
 T_s\in[R_s\cap[h,\infty)]^3
\]
for each \(s\subseteq u\). There are only finitely many such \(s\).
Choose \(M\) larger than every point occurring in the triples \(T_s\),
and put \(a=F\cap\NatLt M\). If \(x\in T_s\), then
\[
 (s\cup a)\cap\NatLt{x+1}
   =(s\cup F)\cap\NatLt{x+1}.
\]
The prefix condition \eqref{eq:bi-prefix} in the definition of
admissibility therefore gives
\(x\in\theta(V_{u,h}^i,s\cup a)\). Consequently, for this single
finite \(a\subseteq F\), the set
\(\theta(V_{u,h}^i,s\cup a)\cap[h,\infty)\) contains the triple \(T_s\)
for every \(s\subseteq u\). Thus no \(s\subseteq u\) witnesses the
right-hand side of \eqref{eq:rapid-11} for this \(a\), a contradiction.
\end{claimproof}

Form the diagonal union
\begin{equation}\label{eq:rapid-12}
 \begin{split}
 L_i=\{\varnothing\}\cup
 \{F\in[H_i]^{<\infty}\setminus\{\varnothing\}:\
 &\exists h\le\min F\\
 &\exists u\subseteq H\cap\NatLt{h}\quad F\in L_{u,h}^i\}.
\end{split}
\end{equation}
\begin{claim}\label{clm:rapid-Li-compact}
\(L_i\) is a compact hereditary subfamily of \(2^{H_i}\), equipped
with the product topology, and every member of \(L_i\) is finite.
\end{claim}
\begin{claimproof}
Heredity and finiteness follow directly from \eqref{eq:rapid-12} and
the heredity of each \(L_{u,h}^i\). It remains to prove closedness in
\(2^{H_i}\). Suppose that \(F_n\in L_i\) and \(F_n\to F\). If
\(F=\varnothing\), then \(F\in L_i\) by definition. Otherwise, put
\(m=\min F\). For all sufficiently large \(n\), we have
\(\min F_n=m\). For each such \(n\), membership in \(L_i\) is witnessed
by parameters
\[
 h_n\le m,\qquad u_n\subseteq H\cap\NatLt{h_n},\qquad
 F_n\in L_{u_n,h_n}^i.
\]
There are only finitely many possible pairs \((u_n,h_n)\), so some
fixed pair \((u,h)\) occurs along a subsequence. Since \(L_{u,h}^i\)
is closed, the limit \(F\) belongs to \(L_{u,h}^i\).
Claim~\ref{clm:rapid-local-no-infinite} then shows that \(F\) is
finite, and \eqref{eq:rapid-12} gives
\(F\in L_i\). Thus \(L_i\) is closed in the compact space \(2^{H_i}\)
and hence is compact.
\end{claimproof}

Membership of a finite set in \eqref{eq:rapid-12} is computable in
\(W\oplus H_i\). Indeed, fix a finite \(F\subseteq H_i\). If
\(F=\varnothing\), membership is immediate from \eqref{eq:rapid-12},
so suppose that \(F\ne\varnothing\). There are only finitely many
possible pairs \((u,h)\) with \(h\le\min F\) and
\(u\subseteq H\cap\NatLt{h}\). For each such pair,
\(V_{u,h}^i\) is computable in \(W\oplus H_i\), since \(Z\le_TW\) and
\[
 V_{u,h}^i=u\cup(H_i\cap[h,\infty)).
\]
The functional \(\Phi\) fixed in
Subsection~\ref{subsec:rapid-responses} therefore computes, for every
\(a\subseteq F\) and \(s\subseteq u\),
\[
 \Phi^{Z\oplus V_{u,h}^i}(\ulcorner s\cup a\urcorner)
   =\ulcorner\theta(V_{u,h}^i,s\cup a)\urcorner.
\]
This is a code for a finite set. Since \(j\) is \(Z\)-computable, we
can decide whether its part above \(h\) contains a \(j\)-zero triple.
The quantifiers over \(a\subseteq F\), \(s\subseteq u\), and the
finitely many pairs \((u,h)\) are all finite, so \eqref{eq:rapid-11} and then
\eqref{eq:rapid-12} are decidable in \(W\oplus H_i\).

For a compact hereditary family \(K\subseteq[H_i]^{<\infty}\), set
\[
 U(K)=\{F_0\cup F_1:F_0,F_1\in K\}.
\]
Define
\[
 \operatorname{Ram}(K)
 =\bigl\{F\in[H_i]^{<\infty}:
 \exists c\colon[F]^2\to2\ \forall G\subseteq F\,
 (|c``[G]^2|\le1\Longrightarrow G\in K)\bigr\}.
\]
\begin{claim}\label{clm:rapid-closure-operations}
Both \(U\) and \(\operatorname{Ram}\) preserve compactness, heredity,
and the absence of infinite members.
\end{claim}
\begin{claimproof}
The union map
\[
 K\times K\longrightarrow 2^{H_i},\qquad(F_0,F_1)\longmapsto F_0\cup F_1,
\]
is continuous, so \(U(K)\) is compact. If
\(E\subseteq F_0\cup F_1\), then
\[
 E=(E\cap F_0)\cup(E\cap F_1),
\]
and both intersections belong to \(K\) by heredity. Thus \(U(K)\) is
hereditary. Its members are finite unions of finite sets, so it has no
infinite member.

For \(\operatorname{Ram}(K)\), work in the compact space
\[
 2^{H_i}\times2^{[H_i]^2}
\]
and let \(\mathcal C_K\) consist of the pairs \((F,c)\) such that, for
every finite \(G\subseteq F\),
\[
 |c([G]^2)|\le1\quad\Longrightarrow\quad G\in K.
\]
This set is closed. Indeed, failure is witnessed by a finite
\(G\notin K\) such that \(G\subseteq F\) and \(c\) is homogeneous on
\([G]^2\), and these two conditions involve only finitely many
coordinates of \(F\) and \(c\). Hence the projection
\[
 P_K=\{F\subseteq H_i:\exists c\in2^{[H_i]^2}\ (F,c)\in\mathcal C_K\}
\]
is compact.

We claim that \(P_K\) has no infinite member. Otherwise, let
\((F,c)\in\mathcal C_K\) with \(F\) infinite. By the infinite Ramsey
theorem, \(F\) has an infinite \(c\)-homogeneous subset \(Y\). Every
finite subset of \(Y\) then belongs to \(K\). An increasing sequence
of finite initial segments of \(Y\) converges to \(Y\) in \(2^{H_i}\),
so the closedness of \(K\) gives \(Y\in K\), contrary to the assumption
that \(K\) has no infinite member.

For finite \(F\), a coloring of \([F]^2\) extends arbitrarily to a
coloring of \([H_i]^2\). It follows that
\[
 P_K=\operatorname{Ram}(K).
\]
Thus \(\operatorname{Ram}(K)\) is compact and has no infinite member.
Finally, restricting a witnessing coloring from \(F\) to a subset of
\(F\) proves heredity.
\end{claimproof}

Starting with \(K_{i,0}=L_i\cup[H_i]^{\le1}\), set
\begin{equation}\label{eq:rapid-13}
 \begin{split}
 K_{i,n+1}
 &=K_{i,n}\cup U(K_{i,n})\cup\operatorname{Ram}(K_{i,n}),\\
 \nu_i(F)&=\min\{n:F\in K_{i,n}\}.
 \end{split}
\end{equation}
All finite sets occur, and membership is computable in
\(W\oplus H_i\). Moreover,
\begin{equation}\label{eq:rapid-14}
 \nu_i(F\cup F')\le\max(\nu_i(F),\nu_i(F'))+1.
\end{equation}
If \(\nu_i(F)>n+1\), every binary pair coloring of \(F\) has a
homogeneous subset of norm greater than \(n\).

By \eqref{eq:rapid-12}, for every \(h\in\N\) and
\(u\subseteq H\cap\NatLt{h}\),
\[
 L_{u,h}^i\subseteq L_i\subseteq K_{i,0}.
\]
Consequently, \(\nu_i(F)=0\) for every \(F\in L_{u,h}^i\).

\subsection{First response $A$}\label{subsec:rapid-freezing}
Using the virtual norms, we build the first homogeneous response \(A\) in
blocks and freeze the relevant finite outputs in order to localize all
later outputs.

\begin{claim}\label{clm:rapid-large-virtual-block}
For every \(i\in I\) and \(n,h\in\N\), there is an \(N\ge h\) such
that
\[
 \nu_i(H_i\cap[h,N])>n.
\]
Moreover, the least such \(N\) is computable in \(W\oplus H_i\),
uniformly in \(n\) and \(h\).
\end{claim}
\begin{claimproof}
If no such \(N\) existed, then
\(\nu_i(H_i\cap[h,N])\le n\), and hence
\(H_i\cap[h,N]\in K_{i,n}\), for every \(N\ge h\). These finite sets
converge in \(2^{H_i}\) to the infinite
set \(H_i\cap[h,\infty)\). Since \(K_{i,n}\) is compact and hence
closed in \(2^{H_i}\), the limit would belong to \(K_{i,n}\),
contrary to the fact that \(K_{i,n}\) has no infinite member. Thus the
required \(N\) exists. Membership in \(K_{i,n}\) is computable in
\(W\oplus H_i\), so searching for the least \(N\) terminates and gives
the asserted computation.
\end{claimproof}

Fix a schedule \(\iota\colon\N\to I\) in which every label $i \in I$ occurs
infinitely often.
We recursively construct finite sets
\[
 (u_m)_{m\in\N},\qquad (P_m)_{m\in\N},\qquad (Q_m)_{m\in\N},
\]
and numbers
\[
 (h_m)_{m\in\N},\qquad (N_m)_{m\in\N},\qquad (M_m)_{m\in\N}.
\]
Here \(u_m\) is the stem available at stage \(m\), with
\(u_0=\varnothing\) and, for \(m>0\),
\(
 u_m=\bigcup_{\ell<m}Q_\ell.
\)
The set \(P_m\) is the allowed second-input block, \(Q_m\) is the
corresponding output block, and \(h_m\le N_m\le M_m\) are their cutoffs.
Put \(M_{-1}=-1\). At stage
\(m\), suppose that the finite \(r\)-\(\varepsilon\)-homogeneous stem
\(u_m\) has been defined, put \(i=\iota(m)\), and let
\[
 n_m=|\{\ell<m:\iota(\ell)=i\}|.
\]
Thus stage \(m\) is the \(n_m\)-th occurrence of label \(i\), counting
from zero. For each proposed starting point \(h>M_{m-1}\), use the
oracle \(W\oplus H_i\) as follows.

\begin{enumerate}
\item
  Let \(N_m(h)\ge h\) be the least number such that
  \begin{equation}\label{eq:rapid-15}
   P_m(h)=H_i\cap[h,N_m(h)],\qquad \nu_i(P_m(h))>n_m.
  \end{equation}
  Claim~\ref{clm:rapid-large-virtual-block} guarantees that \(N_m(h)\)
  is defined.
\item
  For every \(b\subseteq u_m\cup P_m(h)\), compute the entire finite
  output \(\theta(V_{u_m,h}^i,b)\), together with a prefix of
  \(V_{u_m,h}^i\) fixing that output on \([H]^\N\).
\item
  Choose \(M_m(h)\ge N_m(h)\) beyond all the output coordinates and all
  the prefix bounds in step 2.
\end{enumerate}

This defines a \(W\oplus H_i\)-computable function
\(h\mapsto M_m(h)\) for \(h>M_{m-1}\), with \(M_m(h)\ge h\). Extend it
to the finitely many remaining arguments by setting \(M_m(h)=h\).

For every \(x\in u_m\), there is an \(\ell<m\) such that
\(x\in Q_\ell\subseteq H_{\iota(\ell)}\). Since
\(\iota(\ell)\in I\), the \(r\)-row of \(x\) has limit
\(\varepsilon\) along \(X\). As \(u_m\) is finite and \(H_i\subseteq X\),
there is a \(t_m>M_{m-1}\) such that, for every \(h\ge t_m\),
\begin{equation}\label{eq:rapid-cross-edges}
 r(x,y)=\varepsilon
 \quad\text{for every \(x\in u_m\) and
 \(y\in H_i\cap[h,\infty)\)}.
\end{equation}
Apply \eqref{eq:rapid-9} and the fourth fact in
Subsection~\ref{subsec:rapid-candidates} to the extended function
\(h\mapsto M_m(h)\), with lower bound \(t_m\). This gives a
\(h\ge t_m\) such that some \(e\le h\) enters \(J\) after \(M_m(h)\).

Fix such a \(h\) as \(h_m\), and put
\begin{equation}\label{eq:rapid-16}
 \begin{split}
 N_m&=N_m(h_m),\qquad P_m=P_m(h_m),\\
 M_m&=M_m(h_m),\qquad Q_m=H_i\cap[h_m,M_m].
 \end{split}
\end{equation}
Set \(u_{m+1}=u_m\cup Q_m\). Mark \(P_m\) as an allowed second-input
block and \(Q_m\) as its output block. This completes stage \(m\);
stage \(m+1\) starts above \(M_m\).

\begin{figure}[ht]
\centering
\begin{tikzpicture}[x=.9cm,y=.6cm,>=Latex]
\foreach \a/\b/\c/\m in
  {.4/1.6/2.8/0,3.6/4.9/6.1/1,7/8/9.3/2,10.2/11.3/12.6/3}{
  \fill[gray!15] (\a,-.25) rectangle (\c,.45);
  \fill[gray!50] (\a,-.25) rectangle (\b,.45);
  \node[above,font=\small] at (\b,.5) {$Q_{\m}$};
  \node[below,font=\small] at (\a,-.3) {$P_{\m}$};
}
\draw[->,semithick] (0,0)--(13.4,0);
\end{tikzpicture}
\caption{The block arrangement in the construction of $A$. The darker
part of each $Q_m$ is the second-input block $P_m$.}
\label{fig:rapid-block-arrangement}
\end{figure}
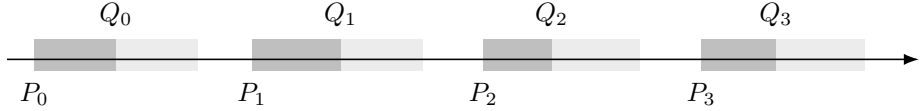

After completing all stages, let
\[
 A=\bigcup_{m\in\N}Q_m=\bigcup_{m\in\N}u_m.
\]
Then
\(A\subseteq H\) is infinite and exactly
\(r\)-\(\varepsilon\)-homogeneous. Indeed, each \(Q_m\) is contained in
an \(r\)-\(\varepsilon\)-homogeneous set \(H_{\iota(m)}\), while
\eqref{eq:rapid-cross-edges} and \(h_m\ge t_m\) make every edge from
\(u_m\) to \(Q_m\) have \(r\)-color \(\varepsilon\).

Let
\begin{equation}\label{eq:rapid-17}
 G_i=\bigcup_{m:\iota(m)=i}P_m,\qquad
 G=\bigcup_{i\in I}G_i\subseteq A.
\end{equation}
For every \(m\) and every \(b\subseteq u_m\cup P_m\), the final first
input \(A\) has the same finite output as the virtual first input:
\begin{equation}\label{eq:rapid-frozen-output}
 \theta(A,b)=\theta(V_{u_m,h_m}^{\iota(m)},b).
\end{equation}
Indeed, let \(\ell_{m,b}\) be the prefix bound computed for \(b\) in
step~2, so that, with \(V=V_{u_m,h_m}^{\iota(m)}\),
\[
 E\cap\NatLt{\ell_{m,b}}=V\cap\NatLt{\ell_{m,b}}
 \quad\Longrightarrow\quad
 \theta(E,b)=\theta(V,b)
 \qquad(E\in[H]^\N).
\]
By the choice of \(M_m\), we have \(\ell_{m,b}\le M_m+1\) and
\(\theta(V,b)\subseteq\NatLt{M_m}\). Since every block added after
stage \(m\) starts above \(M_m\),
\[
 A\cap\NatLt{M_m+1}
 =u_m\cup\bigl(H_{\iota(m)}\cap[h_m,M_m]\bigr)
 =V\cap\NatLt{M_m+1}.
\]
Thus the preceding implication with \(E=A\) proves
\eqref{eq:rapid-frozen-output}, and also
\(\theta(A,b)\subseteq\NatLt{M_m}\).

\begin{claim}[Localization]\label{clm:rapid-localization}
If
\(B\subseteq G\), \(x\in\theta(A,B)\), and \(P_m\) is the last
input block met by \(B\cap\NatLt{x+1}\), then
\begin{equation}\label{eq:rapid-18}
 x\in Q_m.
\end{equation}
\end{claim}
\begin{claimproof}
Put \(b=B\cap\NatLt{x+1}\). Since
\(x\in\theta(A,B)\cap\NatLt{x+1}\), the left-hand side of
\eqref{eq:bi-count}, with \(n=x\), is positive. Hence
\(b\ne\varnothing\). Moreover, applying the prefix condition
\eqref{eq:bi-prefix} first to
\((A,B)\) and then to \((A,b)\) gives
\[
 x\in\theta(A,B)
 \quad\Longleftrightarrow\quad
 x\in\theta(A\cap\NatLt{x+1},b)
 \quad\Longleftrightarrow\quad
 x\in\theta(A,b).
\]

By assumption, \(b\) meets \(P_m\) and no block \(P_\ell\) with
\(\ell>m\). Every earlier block satisfies
\(P_\ell\subseteq Q_\ell\subseteq u_m\), so
\[
 b\subseteq u_m\cup P_m,
 \qquad b\cap P_m\ne\varnothing.
\]
Thus \(b\) is one of the finite inputs frozen at stage \(m\).
Equation~\eqref{eq:rapid-frozen-output} and the output bound following
it yield
\[
 x\in\theta(A,b)
   =\theta(V_{u_m,h_m}^{\iota(m)},b)
   \subseteq\NatLt{M_m},
\]
and hence \(x<M_m\). On the other hand, choose
\(y\in b\cap P_m\). Since \(b\subseteq\NatLt{x+1}\) and
\(P_m\subseteq[h_m,M_m]\), we have
\[
 h_m\le y\le x.
\]
Finally, the support condition \eqref{eq:bi-support} gives
\(x\in A\cup B=A\), because
\(B\subseteq G\subseteq A\). Therefore
\[
 x\in A\cap[h_m,M_m]
   =H_{\iota(m)}\cap[h_m,M_m]
   =Q_m,
\]
as required.
\end{claimproof}

Thus a second input eventually contained in \(G_i\) has almost all its
output in \(H_i\). Moreover, by \eqref{eq:rapid-5},
\begin{equation}\label{eq:rapid-19}
 j_\infty\restriction[Q_m]^2\equiv1.
\end{equation}
Equation~\eqref{eq:rapid-19} controls only the limiting pair colors
and does not imply that all triples in \(Q_m\) have \(j\)-color zero.
Later, the norm arguments choose finite second-input pieces whose
outputs contain \(j\)-zero triples localized in the corresponding
blocks \(Q_m\).

\subsection{Actual norms $(\nu^*_i)_i$}\label{subsec:rapid-transfer}
Once the first response \(A\) is fixed, we transfer the virtual norms
to its actual finite-response map and show that the two norm systems
determine the same block ideal.

Retain from Subsection~\ref{subsec:rapid-freezing} the schedule
\(\iota\), the finite stems \(u_m\), the endpoints \(h_m,N_m,M_m\),
the blocks \(P_m,Q_m\), and the resulting sets \(A,(G_i)_{i\in I},G\).
We make no computability assumption on the finite-response map
\(b\mapsto\theta(A,b)\) in a single row oracle.

For \(i\in I\), a nonempty finite family
\(\mathcal S\subseteq[G]^{<\infty}\), and \(t\in\N\), define a
family on \(G_i\) by
\begin{equation}\label{eq:rapid-20}
 \begin{split}
 D^i_{\mathcal S,t}=\{F\subseteq G_i:
 \forall a\in[F]^{<\infty}\ \exists s\in\mathcal S,\quad&
 \theta(A,s\cup a)\cap g^{-1}(i)\cap[t,\infty)\\
 &\text{contains no \(j\)-zero triple}\}.
 \end{split}
\end{equation}
This family is hereditary and compact.
\begin{claim}\label{clm:rapid-D-no-infinite}
\(D^i_{\mathcal S,t}\) has no infinite member.
\end{claim}
\begin{claimproof}
Suppose that \(F\in D^i_{\mathcal S,t}\) were infinite, and for
\(s\in\mathcal S\) put
\[
 R_s=\theta(A,s\cup F).
\]
Admissibility makes every \(R_s\) infinite. Let
\[
 \ell_s=\max\{m:s\cap P_m\ne\varnothing\},
 \qquad(\max\varnothing=-1).
\]
Since \(F\subseteq G_i\) is infinite and every \(P_m\) is finite, we
may choose \(y_s\in F\cap P_{m_s}\) with \(m_s>\ell_s\). For
\(x\in R_s\cap[y_s,\infty)\), the last block met by
\((s\cup F)\cap\NatLt{x+1}\) is therefore some \(P_m\subseteq G_i\).
Localization \eqref{eq:rapid-18} then yields
\[
 R_s\cap[y_s,\infty)\subseteq g^{-1}(i),
 \qquad\text{and hence}\qquad
 R_s\setminus g^{-1}(i)\text{ is finite}.
\]
Consequently \(R_s\cap g^{-1}(i)\cap[t,\infty)\) is infinite. As
\(j\) has no infinite homogeneous set of color \(1\), choose a
\(j\)-zero triple
\[
 T_s\in[R_s\cap g^{-1}(i)\cap[t,\infty)]^3
 \qquad(s\in\mathcal S).
\]
Choose \(M>\max\bigcup_{s\in\mathcal S}T_s\) and set
\(a=F\cap\NatLt M\). For every \(s\in\mathcal S\) and \(x\in T_s\),
\[
 (s\cup a)\cap\NatLt{x+1}
 =(s\cup F)\cap\NatLt{x+1},
\]
so the prefix condition \eqref{eq:bi-prefix} implies
\[
 T_s\subseteq
 \theta(A,s\cup a)\cap g^{-1}(i)\cap[t,\infty).
\]
Thus this single \(a\in[F]^{<\infty}\) admits no witness
\(s\in\mathcal S\) in \eqref{eq:rapid-20}, a contradiction.
\end{claimproof}

For a family \(\mathcal F\subseteq2^{H_i}\) and a set
\(P\subseteq H_i\), write
\[
 \mathcal F\restriction P=\{F\subseteq P:F\in\mathcal F\}.
\]

\begin{claim}\label{clm:rapid-key-transfer}
For every \(i\in I\), every nonempty finite
\(\mathcal S\subseteq[G]^{<\infty}\), and every \(t\in\N\),
\begin{equation}\label{eq:rapid-21}
 D^i_{\mathcal S,t}\restriction P_m
 \ \subseteq\ K_{i,0}\restriction P_m
 \quad\text{for all sufficiently late \(m\) with \(\iota(m)=i\)}.
\end{equation}
\end{claim}
\begin{claimproof}
For all sufficiently large \(m\) with \(\iota(m)=i\), the first stem
\(u_m\) contains \(\bigcup\mathcal S\), and \(h_m>t\). If
\(F\subseteq P_m\), then \eqref{eq:rapid-frozen-output} gives
\[
 \theta(A,s\cup a)
 =\theta(V^i_{u_m,h_m},s\cup a)
\]
for every \(s\in\mathcal S\) and \(a\subseteq F\). By the support
condition \eqref{eq:bi-support}, all points of the virtual output at or
above \(h_m\) belong to row \(i\). If \(F\in D^i_{\mathcal S,t}\),
then each \(a\subseteq F\) has a witness \(s\in\mathcal S\subseteq
\mathcal P(u_m)\) for which those new outputs contain no \(j\)-zero
triple. Hence
\[
 F\in L^i_{u_m,h_m}; \text{ therefore, }F \in  L_i\subseteq K_{i,0},
\]
which proves \eqref{eq:rapid-21}.
\end{claimproof}

For each \(i\in I\), enumerate all the families
\(D^i_{\mathcal S,t}\) from \eqref{eq:rapid-20} as
\(D_{i,0},D_{i,1},\ldots\), and define compact hereditary families on
\(G_i\) by
\[
 K^*_{i,0}=K_{i,0}\restriction G_i,
\]
\begin{equation}\label{eq:rapid-22a}
 \begin{split}
 K^*_{i,n+1}
 ={}&K^*_{i,n}\cup(K_{i,n+1}\restriction G_i)\cup D_{i,n}\\
 &{}\cup U(K^*_{i,n})\cup\operatorname{Ram}(K^*_{i,n}).
\end{split}
\end{equation}
Define the resulting norm on \([G_i]^{<\infty}\) by
\begin{equation}\label{eq:rapid-22b}
 \nu_i^*(F)=\min\{n:F\in K^*_{i,n}\}.
\end{equation}
The union and Ramsey norm inequalities hold with an increase of one
level, exactly as in \eqref{eq:rapid-14}.

For each \(i\in I\), define two block ideals on \(G_i\) by
\[
 \begin{aligned}
 \mathcal I_i
 &=\left\{Y\subseteq G_i:
   \sup_{m:\iota(m)=i}\nu_i(Y\cap P_m)<\infty\right\},\\
 \mathcal I_i^*
 &=\left\{Y\subseteq G_i:
   \sup_{m:\iota(m)=i}\nu_i^*(Y\cap P_m)<\infty\right\}.
\end{aligned}
\]
The union inequalities show that these are indeed ideals.
\begin{claim}\label{clm:rapid-block-ideals}
For every \(i\in I\),
\begin{equation}\label{eq:rapid-22}
 \mathcal I_i=\mathcal I_i^*.
\end{equation}
\end{claim}
\begin{claimproof}
First, \eqref{eq:rapid-22a} gives, for every \(n\),
\[
 K_{i,n}\restriction G_i\subseteq K^*_{i,n}.
\]
In the other direction, we show by induction on \(n\) that
\begin{equation}\label{eq:rapid-level-transfer}
 K^*_{i,n}\restriction P_m
 \subseteq K_{i,n}\restriction P_m
 \quad\text{for all sufficiently late \(m\) with \(\iota(m)=i\)}.
\end{equation}
For \(n=0\), this follows from
\(K^*_{i,0}=K_{i,0}\restriction G_i\). For the induction step, take
\(m\) late enough that both the induction hypothesis and
Claim~\ref{clm:rapid-key-transfer} apply. The term
\(K_{i,n+1}\restriction G_i\) in \eqref{eq:rapid-22a} already lies at
virtual level \(n+1\), and the new generator \(D_{i,n}\) lies on
\(P_m\) in \(K_{i,0}\) by Claim~\ref{clm:rapid-key-transfer}. For the
old term \(K^*_{i,n}\), the induction hypothesis and
\(K_{i,n}\subseteq K_{i,n+1}\) apply. Finally, the induction hypothesis
applies to every subset of \(P_m\); the definitions of \(U\) and
\(\operatorname{Ram}\), together with \eqref{eq:rapid-13}, then put
their restrictions to \(P_m\) in \(K_{i,n+1}\). This proves
\eqref{eq:rapid-level-transfer}.

Now suppose \(Y\in\mathcal I_i\), and choose \(N\) such that
\[
 \nu_i(Y\cap P_m)\le N
 \quad\text{whenever \(\iota(m)=i\)}.
\]
Then \(Y\cap P_m\in K_{i,N}\restriction G_i\subseteq K^*_{i,N}\) for
every such \(m\), so \(Y\in\mathcal I_i^*\).

Conversely, suppose \(Y\in\mathcal I_i^*\), and choose \(N\) bounding
all the values \(\nu_i^*(Y\cap P_m)\) with \(\iota(m)=i\).
By \eqref{eq:rapid-level-transfer},
\(\nu_i(Y\cap P_m)\le N\) on all sufficiently late row-\(i\) blocks.
There are only finitely many earlier blocks, and each
\(\nu_i(Y\cap P_m)\) is finite. Taking the maximum of \(N\) and these
finitely many values gives a bound on all the virtual norms. Hence
\(Y\in\mathcal I_i\), proving \eqref{eq:rapid-22}.
\end{claimproof}

This common ideal is proper by \eqref{eq:rapid-15}. Call
\(Y\subseteq G_i\) \emph{positive} if \(Y\notin\mathcal I_i\).

The actual norms \(\nu_i^*\) need not be effective. They are used only
after \(c\) and \(A\) have been constructed and an arbitrary second
coloring \(d\in\Ccal^2_2\) has been given. The endpoint functions
\(h\mapsto M_m(h)\) used in the construction of \(A\) depend only on
the effective virtual norms \(\nu_i\).

\subsection{Arbitrary second colorings $d$}\label{subsec:rapid-second}
We use the actual norms to handle any prescribed \(d\in\Ccal^2_2\).
We construct an infinite exactly \(d\)-homogeneous \(B\subseteq G\)
for which \(\theta(A,B)\) is not \(c\)-almost homogeneous. We
distinguish whether some \(G_i\) already contains a positive
\(d\)-homogeneous set.
Claims~\ref{clm:rapid-distributed-case} and
\ref{clm:rapid-single-Gi-case} treat the two alternatives: first the
case with no such \(i\), using finite pieces from infinitely many
distinct sets \(G_i\), and then the case in which one \(G_i\) suffices.

\begin{claim}\label{clm:rapid-color-switch}
Suppose \(Y\subseteq G_i\) is positive and all finite
\(d\)-homogeneous subsets of \(Y\) of color \(0\) belong to
\(K^*_{i,N}\). Then \(Y\) contains a positive \(d\)-homogeneous set
of color \(1\). The same statement holds with the colors reversed.
\end{claim}
\begin{claimproof}
Put \(Y_0=Y\). As long as possible, choose \(x_r\in Y_r\) such that
the forward 0-neighborhood
\[
 Y_{r+1}
 =\{y\in Y_r:y>x_r\text{ and }d(x_r,y)=0\}
\]
is positive. If this continued for every \(r\), then
\(\{x_r:r\in\N\}\) would be an infinite \(d\)-homogeneous set of
color \(0\). Every finite initial segment would belong to
\(K^*_{i,N}\) by hypothesis, so the closedness of \(K^*_{i,N}\)
would put the infinite limit in \(K^*_{i,N}\), contradicting the
absence of infinite members.
Therefore the process stops with a positive set \(Z\subseteq Y\) such
that, for every \(x\in Z\),
\begin{equation}\label{eq:rapid-small-zero-neighborhood}
 Z_x^0=\{y\in Z:y>x\text{ and }d(x,y)=0\}
 \in\mathcal I_i^*.
\end{equation}

We recursively construct nonempty finite sets \(F_m\subseteq Z\).
Suppose that \(F_0,\ldots,F_{m-1}\) have been chosen, and put
\[
 S_m=\bigcup_{\ell<m}F_\ell,\qquad
 Z_m=\{y\in Z:y>\max S_m\text{ and }
              (\forall x\in S_m)\ d(x,y)=1\},
\]
where \(Z_0=Z\). By \eqref{eq:rapid-small-zero-neighborhood}, the set
removed from a tail of \(Z\) is contained in
\(\bigcup_{x\in S_m}Z_x^0\), a finite union of members of
\(\mathcal I_i^*\). Hence \(Z_m\) is positive. Its actual block norms
are therefore unbounded, so, after avoiding the finitely many blocks
used earlier, choose \(P_{\ell_m}\), with \(\iota(\ell_m)=i\), such that
\[
 \nu_i^*(Z_m\cap P_{\ell_m})>\max(N,m)+1.
\]
The Ramsey norm inequality supplies a \(d\)-homogeneous
\(F_m\subseteq Z_m\cap P_{\ell_m}\) satisfying
\[
 \nu_i^*(F_m)>\max(N,m).
\]
Its \(d\)-color cannot be \(0\), since the hypothesis would then give
\(F_m\in K^*_{i,N}\) and hence \(\nu_i^*(F_m)\le N\). Thus \(F_m\)
has \(d\)-color \(1\). The definition of \(Z_m\) also makes every
edge from \(S_m\) to \(F_m\) have \(d\)-color \(1\).

It follows that
\[
 Y'=\bigcup_{m\in\N}F_m
\]
is \(d\)-homogeneous of color \(1\). Since the families \(K^*_{i,n}\)
are hereditary,
\[
 \nu_i^*(Y'\cap P_{\ell_m})
 \ge \nu_i^*(F_m)>m.
\]
Thus \(Y'\notin\mathcal I_i^*\), so \(Y'\) is positive by
Claim~\ref{clm:rapid-block-ideals}. The argument with the colors
reversed is identical.
\end{claimproof}

\begin{claim}\label{clm:rapid-finite-extension}
Assume that no \(G_j\) contains a positive \(d\)-homogeneous set. Let
\(Y\subseteq G_i\) be positive, and suppose that positive
sets \(Y_j\subseteq G_j\) are given for infinitely many
\(j\in I\setminus\{i\}\). For every \(N\), there are a color
\(\delta\in\Bits\), a finite \(d\)-homogeneous set \(F\subseteq Y\) of
color \(\delta\) with \(\nu_i^*(F)>N\), and positive
candidate sets \(Y'_j\subseteq Y_j\) for infinitely many \(j\), such
that every edge from \(F\) to every \(Y'_j\) has color \(\delta\).
\end{claim}
\begin{claimproof}
Claim~\ref{clm:rapid-color-switch} implies that every positive subset
of every \(G_j\) has finite \(d\)-homogeneous subsets of arbitrarily
large actual norm in both colors. Indeed, if the norms in one color
were bounded on a positive set, that claim would produce a positive
homogeneous set in the other color, contrary to the present
hypothesis.

Let \(J_{-1}\subseteq I\setminus\{i\}\) be the infinite set of indices
for which the candidates are given, put
\[
 Z_{-1,j}=Y_j\quad(j\in J_{-1}),
 \qquad Y=\{x_0<x_1<\cdots\}.
\]
Recursively in \(t\), suppose that the infinite set \(J_{t-1}\) and
positive sets \(Z_{t-1,j}\subseteq Y_j\), \(j\in J_{t-1}\), have been
defined. For \(j\in J_{t-1}\) and \(\varepsilon\in\Bits\), set
\[
 Z^{\varepsilon}_{t,j}
 =\{y\in Z_{t-1,j}:d(x_t,y)=\varepsilon\}.
\]
Since
\[
 Z_{t-1,j}=Z^0_{t,j}\cup Z^1_{t,j}
\]
is positive and \(\mathcal I_j^*\) is an ideal, choose
\(\varepsilon_{t,j}\in\Bits\) so that
\(Z^{\varepsilon_{t,j}}_{t,j}\) is positive. By the infinite
pigeonhole principle, there are \(\delta_t\in\Bits\) and an infinite
set
\[
 J_t\subseteq
 \{j\in J_{t-1}:\varepsilon_{t,j}=\delta_t\}.
\]
For \(j\in J_t\), put
\[
 Z_{t,j}=Z^{\delta_t}_{t,j}.
\]
Thus every \(Z_{t,j}\) is positive, and for all \(s\le t\),
\begin{equation}\label{eq:rapid-finite-extension-compatibility}
 j\in J_t\ \text{ and }\ y\in Z_{t,j}
 \quad\Longrightarrow\quad d(x_s,y)=\delta_s.
\end{equation}

For \(\varepsilon\in\Bits\), let
\[
 X_\varepsilon=\{x_t:\delta_t=\varepsilon\}.
\]
Since \(Y=X_0\cup X_1\) is positive, one of these two sets, say
\(X_\delta\), is positive. By the first paragraph, there is a finite
\(d\)-homogeneous set \(F\subseteq X_\delta\) of color \(\delta\)
such that \(\nu_i^*(F)>N\). Put
\[
 T=\max\{t:x_t\in F\},\qquad
 Y'_j=Z_{T,j}\quad(j\in J_T).
\]
The set \(J_T\) is infinite and every \(Y'_j\) is positive. Moreover,
if \(x_t\in F\), \(j\in J_T\), and \(y\in Y'_j\), then
\(t\le T\) and \(\delta_t=\delta\), so
\eqref{eq:rapid-finite-extension-compatibility} gives
\[
 d(x_t,y)=\delta.
\]
Hence every edge from \(F\) to every \(Y'_j\), \(j\in J_T\), has
color \(\delta\), as required.
\end{claimproof}

\begin{claim}\label{clm:rapid-distributed-case}
If no \(G_i\) contains a positive \(d\)-homogeneous set, then there is
an infinite exactly \(d\)-homogeneous set \(B\subseteq G\) such that
\(\theta(A,B)\) is not \(c\)-almost homogeneous.
\end{claim}
\begin{claimproof}
\emph{Construction.}
Each \(G_i\) is positive, and \(I\) is infinite. We construct finite
second-input pieces \(a_n\), distinct indices \(i_n\), colors
\(\delta_n\), infinite sets \(J_n\subseteq I\), and positive candidate
sets \(Y_{n,j}\subseteq G_j\) for \(j\in J_n\).

Start with
\[
 J_0=I,\qquad Y_{0,j}=G_j\ (j\in J_0),\qquad b_0=\varnothing.
\]
At stage \(n\), assume that \(J_n\), the candidates \(Y_{n,j}\), and
\[
 b_n=\bigcup_{\ell<n}a_\ell
\]
have been defined. Choose \(i_n\in J_n\), and choose \(t_n>n\) above
every point of
\[
 b_n\cup\bigcup_{s\subseteq b_n}\theta(A,s).
\]
For each \(j\in J_n\), put
\[
 \widehat Y_{n,j}=Y_{n,j}\cap[t_n,\infty).
\]
These sets remain positive. Choose \(q_n\) such that
\[
 D^{i_n}_{\mathcal P(b_n),t_n}=D_{i_n,q_n}.
\]
By \eqref{eq:rapid-22a}, this family is contained in
\(K^*_{i_n,q_n+1}\).

Apply Claim~\ref{clm:rapid-finite-extension} to
\(\widehat Y_{n,i_n}\), the candidates
\(\widehat Y_{n,j}\) for \(j\in J_n\setminus\{i_n\}\), and
\(N=q_n+1\). We obtain a color \(\delta_n\), a finite
\(d\)-homogeneous set of color \(\delta_n\),
\[
 F_n\subseteq\widehat Y_{n,i_n},
 \qquad \nu_{i_n}^*(F_n)>q_n+1,
\]
an infinite set \(J_{n+1}\subseteq J_n\setminus\{i_n\}\), and positive
surviving candidates \(Y'_{n,j}\subseteq\widehat Y_{n,j}\),
\(j\in J_{n+1}\), such that every edge from \(F_n\) to every
\(Y'_{n,j}\) has color \(\delta_n\).

Since \(F_n\notin K^*_{i_n,q_n+1}\), we have
\(F_n\notin D^{i_n}_{\mathcal P(b_n),t_n}\). Hence some
\(a_n\subseteq F_n\) satisfies
\begin{equation}\label{eq:rapid-23}
 \forall s\subseteq b_n,\quad
 \theta(A,s\cup a_n)\cap g^{-1}(i_n)\cap[t_n,\infty)
 \text{ contains a \(j\)-zero triple}.
\end{equation}
The choice of \(t_n\) shows that \(a_n\ne\varnothing\): for
\(a_n=\varnothing\), every \(\theta(A,s)\) lies below \(t_n\).

Choose \(v_n>t_n\) above every point of \(F_n\) and every output
coordinate of \(\theta(A,s\cup a_n)\), for every \(s\subseteq b_n\), and set
\[
 Y_{n+1,j}=Y'_{n,j}\cap[v_n,\infty)
 \quad(j\in J_{n+1}),\qquad
 b_{n+1}=b_n\cup a_n.
\]
Removing these finite initial segments preserves positivity and
completes stage \(n\).

\smallskip
\noindent\emph{Verification.}
The indices \(i_n\) are distinct because
\(J_{n+1}\subseteq J_n\setminus\{i_n\}\). Moreover, \(a_n\) is
internally \(d\)-homogeneous of color \(\delta_n\), and, whenever
\(n<m\), the construction of the surviving candidates gives
\[
 d(x,y)=\delta_n
 \quad\text{for every \(x\in a_n\) and \(y\in a_m\)}.
\]
Choose an infinite \(R\subseteq\N\) and a color \(\delta\) such that
\(\delta_n=\delta\) for every \(n\in R\), and put
\[
 B=\bigcup_{n\in R}a_n.
\]
Then \(B\) is infinite and exactly \(d\)-homogeneous of color
\(\delta\).

For \(n\in R\), let
\[
 s_n=\bigcup_{\substack{\ell\in R\\ \ell<n}}a_\ell.
\]
Then \(s_n\subseteq b_n\), so \eqref{eq:rapid-23} supplies a
\(j\)-zero triple \(T_n\) in
\[
 \theta(A,s_n\cup a_n)\cap g^{-1}(i_n)\cap[t_n,\infty).
\]
All pieces \(a_\ell\) with \(\ell>n\) lie above \(v_n\), hence above
every output coordinate of \(\theta(A,s_n\cup a_n)\). The prefix condition
\eqref{eq:bi-prefix} therefore preserves \(T_n\) in
\(\theta(A,B)\).

We now see explicitly that both \(c\)-colors occur arbitrarily far out
in the output. Fix \(L\in\N\). Choose \(n\in R\) with \(n>L\). Since
\[
 T_n\subseteq g^{-1}(i_n)\cap[t_n,\infty)
 \quad\text{and}\quad t_n>n,
\]
the \(c\)-zero triple \(T_n\) lies entirely above \(L\). Choose
\(x\in T_n\), then choose \(m\in R\) with \(m>x\) and \(y\in T_m\).
We have
\[
 x< m<t_m\le y,
 \qquad g(x)=i_n\ne i_m=g(y).
\]
Finally, choose \(\ell\in R\) with \(\ell>y\) and \(z\in T_\ell\). Then
\(x<y<z\), all three points belong to \(\theta(A,B)\), and
\eqref{eq:rapid-6} gives
\[
 c(x,y,z)=1
\]
because \(x\) and \(y\) lie in different fibers. Thus, above every
\(L\), the output contains both the \(c\)-zero triple \(T_n\) and a
\(c\)-one triple \(\{x,y,z\}\). Hence \(\theta(A,B)\) is not
\(c\)-almost homogeneous.
\end{claimproof}

\begin{claim}\label{clm:rapid-single-Gi-case}
If some \(G_i\) contains a positive \(d\)-homogeneous set \(Y\), then
there is an infinite exactly \(d\)-homogeneous set \(B\subseteq Y\)
such that \(\theta(A,B)\) is not \(c\)-almost homogeneous.
\end{claim}
\begin{claimproof}
We recursively construct nonempty finite sets \(a_n\subseteq Y\) and
strictly increasing block indices \(m_n\) with \(\iota(m_n)=i\). Put
\[
 b_0=\varnothing,\qquad b_n=\bigcup_{\ell<n}a_\ell.
\]
At stage \(n\), choose \(t_n>n\) above every point of
\[
 b_n\cup\theta(A,b_n).
\]
Choose \(q_n\) such that
\[
 D^i_{\{b_n\},t_n}=D_{i,q_n}.
\]
By \eqref{eq:rapid-22a},
\(D^i_{\{b_n\},t_n}\subseteq K^*_{i,q_n+1}\).

Since \(Y\) is positive, its actual block norms are unbounded. We may
therefore choose \(m_n>m_{n-1}\), with \(\iota(m_n)=i\) and
\(h_{m_n}>t_n\), such that
\[
 \nu_i^*(Y\cap P_{m_n})>q_n+1,
\]
where the condition \(m_n>m_{n-1}\) is omitted when \(n=0\). It follows
that
\[
 Y\cap P_{m_n}\notin D^i_{\{b_n\},t_n}.
\]
By the definition \eqref{eq:rapid-20}, there is a set
\(a_n\subseteq Y\cap P_{m_n}\) such that
\[
 \theta(A,b_n\cup a_n)\cap g^{-1}(i)\cap[t_n,\infty)
\]
contains a \(j\)-zero triple; fix one and call it \(T_n\). The set
\(a_n\) is nonempty, since \(\theta(A,b_n)\) lies below \(t_n\). Set
\[
 b_{n+1}=b_n\cup a_n.
\]

We next locate \(T_n\). If \(x\in T_n\) and \(x<h_{m_n}\), then
\[
 b_{n+1}\cap\NatLt{x+1}=b_n\cap\NatLt{x+1},
\]
so the prefix condition \eqref{eq:bi-prefix} would give
\(x\in\theta(A,b_n)\), contrary to \(x\ge t_n\). Hence every point of
\(T_n\) is at least \(h_{m_n}\). Since
\(a_\ell\subseteq P_{m_\ell}\subseteq Q_{m_\ell}\subseteq u_{m_n}\)
for every \(\ell<n\), we have \(b_n\subseteq u_{m_n}\). Together with
\(a_n\subseteq P_{m_n}\),
\eqref{eq:rapid-frozen-output} and the choice of \(M_{m_n}\) put every
point of \(T_n\) at most \(M_{m_n}\). The support condition
\eqref{eq:bi-support}, together with
\(b_{n+1}\subseteq Y\subseteq G\subseteq A\), now gives
\[
 T_n\subseteq A\cap[h_{m_n},M_{m_n}]=Q_{m_n}.
\]
Thus \(T_n\) is a \(c\)-zero triple, and
\eqref{eq:rapid-19} gives
\[
 j_\infty\restriction[T_n]^2\equiv1.
\]

Finally, put
\[
 B=\bigcup_{n\in\N}a_n.
\]
The sets \(a_n\) are nonempty and occur successively inside the
\(d\)-homogeneous set \(Y\), so \(B\subseteq Y\) is infinite and
exactly \(d\)-homogeneous. At stage \(n+1\), the choice of
\(t_{n+1}\) places all later inputs above every coordinate of
\(\theta(A,b_{n+1})\).
Therefore \eqref{eq:bi-prefix} preserves \(T_n\) in
\(\theta(A,B)\). Since \(t_n>n\), these give arbitrarily late
\(c\)-zero triples in the output.

For each \(n\), every pair in \(T_n\) has \(j\)-limit \(1\). Taking a
later output point beyond the corresponding stabilization bound gives
a \(c\)-one triple, because the first two points lie in the same fiber
\(g^{-1}(i)\). Hence \(\theta(A,B)\) has arbitrarily late triples of
both colors and is not \(c\)-almost homogeneous.
\end{claimproof}

Claims~\ref{clm:rapid-distributed-case} and
\ref{clm:rapid-single-Gi-case} exhaust the two cases and prove
\eqref{eq:rapid-1}.\hfill\(\square\)

\subsection{Conclusion}

\begin{theorem}\label{thm:two-stage-negative}
For every finite positive integer \(k\), there are no maps
\[
 \alpha_1:\Ccal^3_2\to\Ccal^2_k,\qquad
 \alpha_2:\Ccal^3_2\times\X\to\Ccal^2_2,\qquad
 \Gamma:\X^2\to\X
\]
such that \(\alpha_1,\Gamma\) are Borel and
\[
 \alpha_1(c)H^2_kA\ \land\ \alpha_2(c,A)H^2_2B
 \quad\Longrightarrow\quad cH^3_2\Gamma(A,B)
\]
for all \(c,A,B\). No regularity assumption on \(\alpha_2\) is needed.
In particular,
\[
 \Hom^2_k;\Hom^2_2\not\longrightarrow_{\BT}\Hom^3_2.
\]
\end{theorem}
\begin{proof}
Apply Proposition~\ref{prop:rapid-normal-form} to a proposed reduction.
Theorem~\ref{thm:rapid-impossibility} supplies \(c,A\) and then an
exactly homogeneous response \(B\) to the actual second query, for
which the output is not almost homogeneous. Exact homogeneity implies
almost homogeneity, so both responses are valid and give the required
contradiction.
\end{proof}

\section{Stage count}\label{sec:consequences}

\begin{corollary}[Exactly three binary pair stages]\label{cor:exact-stage-count}
The least positive integer \(k\) for which a sequential composition
of \(k\) copies of \(\Hom^2_2\) Borel--Tukey reduces to
\(\Hom^3_2\) is \(3\).
\end{corollary}
\begin{proof}
Theorem~\ref{thm:three-stage} gives the upper bound, and
Theorem~\ref{thm:two-stage-negative} excludes two stages.
\end{proof}

\section*{Declaration of generative AI and AI-assisted technologies}
During the preparation of this work, the author used OpenAI Codex with
the \texttt{gpt-5.6-sol} and \texttt{gpt-6-astra} models to assist with
literature searches, the discovery and development of proofs, language
editing, and \LaTeX{} drafting. The author independently checked,
reviewed, and revised all AI-assisted material and takes full
responsibility for the content of the article.

\end{document}